%% file: hadwigerramsey.tex
\documentclass{article}
\usepackage[T1]{fontenc}
\usepackage[utf8]{inputenc}
\usepackage{geometry}
\usepackage{amsmath,amsfonts,amssymb,amsthm,thmtools,mathtools}
\usepackage{enumerate}

\usepackage{graphicx}
\usepackage{tikz}
\usetikzlibrary{calc}
\usepackage{comment}
\usepackage[noabbrev,capitalise]{cleveref}
\usepackage{todonotes}

\usepackage[
  backend=bibtex,style=numeric,
  doi=false,url=false,isbn=false,eprint=true,
  giveninits=false,maxbibnames=99
]{biblatex}
\DeclareDelimFormat[bib,biblist]{nametitledelim}{\addcomma\space}

\DeclareFieldFormat[article]{title}{#1}
\DeclareFieldFormat[article]{volume}{\mkbibbold{#1}}
\DeclareFieldFormat[article]{number}{\mkbibparens{#1}}
\DeclareFieldFormat[article]{pages}{#1}
\renewbibmacro*{in:}{%
  \ifentrytype{article}
    {}
    {\bibstring{in}\printunit{\intitlepunct}}}
\renewbibmacro*{volume+number+eid}{%
  \printfield{volume}%
  \setunit*{\addnbthinspace}%
  \printfield{number}%
  \setunit{\bibeidpunct}%
  \printfield{eid}}
\AtEveryBibitem{\ifentrytype{book}{\clearfield{pages}}{}}

\newtheorem{thm}{Theorem}[section]
\newtheorem{lem}[thm]{Lemma}

\newtheorem{prop}[thm]{Proposition}
\newtheorem{cor}[thm]{Corollary}

\newcommand{\eps}{\varepsilon}
 
\newcommand{\mc}[1]{\mathcal{#1}} 
\newcommand{\bb}[1]{\mathbb{#1}}
\newcommand{\brm}[1]{\operatorname{#1}}

\newcommand{\fs}[2]{\left(\frac{#1}{#2}\right)}
\newcommand{\s}[1]{\left(#1\right)}
\newcommand{\set}[1]{\{#1\}}

\usepackage[noabbrev,capitalise]{cleveref}
\renewcommand{\v}{\textup{\textsf{v}}}
\newcommand{\e}{\textup{\textsf{e}}}

\newcommand{\Prob}{\mathbb{P}}
\DeclareMathOperator{\had}{h}

\definecolor{definitionred}{RGB}{128,0,32}

\newcommand{\defn}[1]{\textcolor{definitionred}{\emph{#1}}}

\crefname{thm}{Theorem}{Theorems}
\crefname{lem}{Lemma}{Lemmas}
\crefname{cor}{Corollary}{Corollaries}
\crefname{obs}{Observation}{Observations}

	\title{Asymptotics of the three-color Hadwiger-Ramsey numbers}
	\author{ Sergey Norin\thanks{Department of Mathematics and Statistics, McGill University. Email: {\tt sergey.norin@mcgill.ca}. Supported by an NSERC Discovery grant.} \and Neil Rahman  \thanks{Department of Mathematics and Statistics, McGill University. Email: {\tt neil.rahman@mcgill.ca}. Partially supported by an NSERC Discovery grant and supported by an NSERC CGRS M award.}}

\begin{document}
	\maketitle
	
	\begin{abstract}
		The Hadwiger-Ramsey number $R_h(k;\ell)$ is defined as the smallest integer $n$ such that for any edge coloring of a complete graph on $n$ vertices in $\ell$ colors at least one of the color classes induces a graph with a $K_k$ minor. We  show that $$R_h(k;3) = \s{\beta^{-1} + o(1)}k \sqrt{\ln k},$$
		where $\beta^{-1} = 2/\sqrt{-\ln(10-4\sqrt6)} = 1.58149\ldots$. 
		This resolves the first open case of a problem by Axenovich and Steiner, who determined the asymptotic behavior of $R_h(k;2)$ and asked about the asymptotics of $R_h(k;\ell)$ for $\ell \geq 3$.
		\end{abstract}
	
\section{Introduction}

An \defn{$\ell$-decomposition of the complete graph $K_n$} is a collection of spanning subgraphs $(G_1,\ldots,G_{\ell})$ of $K_n$ such that $(E(G_1),\ldots,E(G_\ell))$ is a partition of $E(K_n)$. Thus $\ell$-decompositions correspond to colorings of edges of $K_n$.
The \defn{classical Ramsey number $R(k;\ell)$} is defined as the smallest integer $n$ such that for any $\ell$-decomposition $(G_1,\ldots,G_\ell)$ of $K_n$ at least one of the graphs  $G_i$ contains a $K_k$ subgraph. Numerous generalizations of classical Ramsey numbers have been extensively investigated...

In particular, in $1974$ Chartrand and Polimeni~\cite{Chartrand_Polimeni_1974} generalized classical Ramsey numbers to  $f$-Ramsey numbers, where for a monotone graph parameter $f$, unbounded on complete graphs, the \defn{$f$-Ramsey number $R_f(k;\ell)$} is defined as the smallest integer $n$ such that for every $\ell$-decomposition $(G_1,\ldots,G_\ell)$ of $K_n$ we have $f(G_i) \geq k$ for some $1 \leq i \leq \ell$.  Thus $R(k;\ell)=R_{\omega}(k;\ell)$, where $\omega(G)$ denotes the clique number. Chartrand and Polimeni exactly determined the $\chi$-Ramsey numbers, where $\chi$ denotes the chromatic number showing that $$R_{\chi}(k;\ell)=1+(k-1)^{\ell}$$ for all pairs of positive integers $k$ and $\ell$.

The $f$-Ramsey numbers have been further investigated for a variety of natural parameters $f$. For example \begin{itemize} \item $R_{\chi'}(k;\ell)$ where $\chi'$ denotes the chromatic index has been completely determined by Lesniak-Foster and Roberts~\cite{LesniakFoster_Roberts_1977},
	\item Matula~\cite{Matula_1983} considered Ramsey theory for graph connectivity showing $$ 2\ell(k-1)+1 \leq R_{\kappa^{\star}}(k;\ell) \leq \frac{10}{3}\ell(k-1)+1$$ for $k,\ell \geq 2$, where $\kappa^{\star}(G)$ denotes the largest integer $k$ such that $G$ contains a $k$-connected subgraph. \footnote{The corresponding natural parameter $\kappa$ equal to graph connectivity is not monotone. The upper bound can be improved to $\frac{19}{6}\ell (k-1) -k +3$, using the bounds of Bernshteyn and Kostochka~\cite{Bernshteyn_Kostochka_2016} on  the maximum number of edges in a graph with no  $k$-connected subgraph.} 
	\item Joret and Wood~\cite{Joret_Wood_2012} proved that $R_{\brm{tw}}(k;2) = 2 k+1$ for $k \geq 2$.
\end{itemize} 

See also~\cite{Burr_Jacobson_Mihok_Semanisin_1999,Burger_vanVuuren_2011} for further work on generalized Ramsey numbers.

In this paper we investigate the $\had$-Ramsey numbers, where 
the \defn{Hadwiger number $\had(G)$} of a graph $G$ is defined as the maximum $k$ such that $G$ has a $K_k$ minor. Hadwiger numbers are the subject of the famous Hadwiger's conjecture~\cite{Hadwiger43} and have been the subject of extensive investigation for decades.

The $\had$-Ramsey numbers, which we refer to as the \defn{Hadwiger-Ramsey numbers}, were surprisingly first systematically investigated only recently by Axenovich and Steiner~\cite{hadramtwo} , who have determined the asymptotic behavior of the two-color Hadwiger-Ramsey numbers, showing that   \begin{equation}\label{e:2color}R_{\had}(k;2) = (1+o(1))k\sqrt{\log_2{k}}.\end{equation} 

They have left the problem of determining the  asymptotics of $R_{\had}(k;\ell)$ for fixed $\ell \geq 3$ open. Our main result resolves the first  case of this problem.

Let $$\beta = \frac{1}{2}\sqrt{-\ln(10-4\sqrt6)}$$
be the constant ubiquitous in our calculations.

\begin{thm}\label{t:main}  
	$$R_{\had}(k;3) = \s{\beta^{-1} + o(1)}k \sqrt{\ln k}.$$
\end{thm}	

The extremal graphs for both two- and three-color Hadwiger-Ramsey numbers arise from random constructions. The result of Axenovich and Steiner implies that uniformly random colorings asymptotically minimize the largest size of a monochromatic complete minor in a $2$-edge-coloring of a complete graph. 

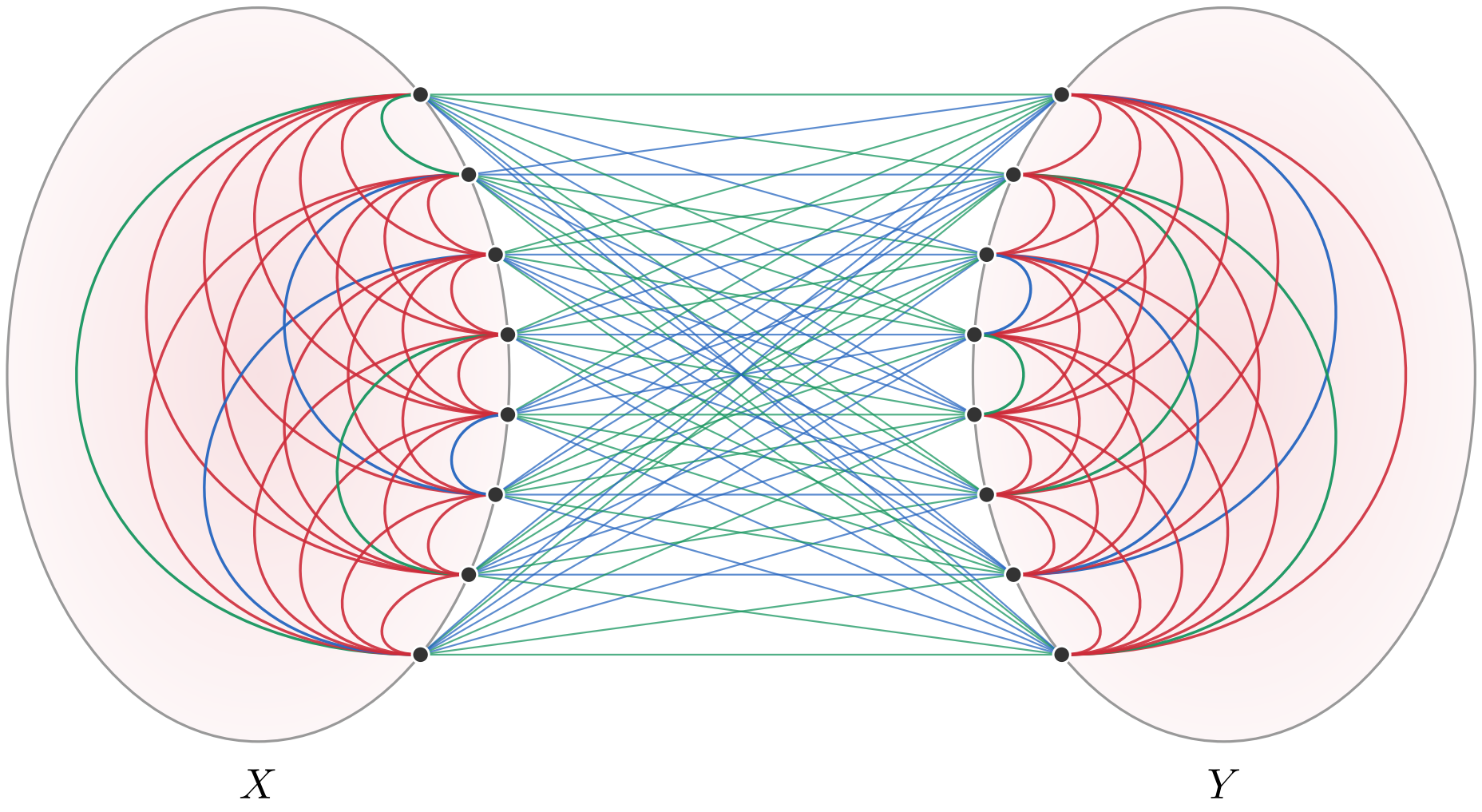
\begin{figure}[t]
	\centering
\input{three-color.tikz}
	\caption{The optimal construction for $3$-color Hadwiger-Ramsey numbers.}
	\label{f:construction}
\end{figure}

One might expect the analogous situation for $3$-edge-coloring, but it turns out that there is a better random construction. We construct $3$-decompositions $(G_1,G_2,G_3)$ of $K_n$ such that $\had(G_i) < k$ for $n=\s{\beta^{-1} + o(1)}k \sqrt{\log k}$, certifying the lower bound on $R_h(k;3)$ in \cref{t:main}, as follows. 

Let $(X,Y)$ be a partition of the vertex set of $K_n$ into two almost equal parts and let $p=5-2\sqrt{6}$. Every edge with both ends in $X$ or both ends in $Y$ is assigned to $G_1$ with probability $1-2p$ and to each of $G_2$ and $G_3$ with probability $p$, independently. Every edge with one end in $X$ and the other in $Y$ is assigned to $G_2$ with probability $1/2$ and, otherwise, to $G_3$. See \cref{f:construction}.

 It immediately follows from our results, although we don't make it explicit, that every $3$-decomposition that is close to optimal for $3$-color Hadwiger-Ramsey numbers is close to the one described above. In particular, the proportion of edges in one of the graphs in the decomposition is close to $\frac{1-2p}{2} \sim 0.399 \neq 1/3$. Thus the parts in optimal decompositions are not symmetric. This is not  unique among known Ramsey constructions, but somewhat unusual.

\subsection{Tools and the proof outline}

The upper and lower bounds in \cref{t:main} correspond to the following two theorems.

\begin{thm}\label{existencethm} For every $\eps > 0$ there exists $n_0$ such that for every  integer $n \geq n_0$  there exists a  $3$-decomposition $(G_1, G_2,G_3)$  of $K_n$ such that 
\begin{equation}\label{e:lower} \had(G_i) \leq (1+ \eps) \beta \frac{n}{\sqrt{\log n}}, \end{equation}  for every $i \in \{1,2,3\}$.
\end{thm}

\begin{thm}\label{mainthm} For every $\eps > 0$ there exists $n_0$ such that for every  integer $n \geq n_0$ and  every $3$-decomposition $(G_1, G_2,G_3)$  of $K_n$ there exists $i \in \{1,2,3\}$  such that $$ \had(G_i) \geq (1-\eps)\beta\frac{n}{\sqrt{\log n}}. $$
\end{thm}

It remains to prove Theorems ~\ref{existencethm} and~\ref{mainthm}. 

Theorem~\ref{existencethm}  is easier. A suitable random $3$-decomposition was described in the introduction and we just need to verify that \eqref{e:lower} holds for every part of it with high probability. The techniques for estimating Hadwiger number of random graphs were first developed by Bollob\'as, Catlin and Erd\H{o}s~\cite{MR593989}. While their result doesn't directly apply in our setting, the techniques extend to the following technical theorem which applies.

\begin{thm}\label{thomasonlower} For every $\eps > 0$ there exists $n_0=n_{0,\ref{thomasonlower}}(\eps)$ such that for all $n \geq n_0$ the following holds. Let $(X_1, \ldots, X_k)$ be a partition of $[n]$, let $0 \leq p_{ij}=p_{ji} \leq 1 -\eps$ for all $1 \leq i,j \leq k$. Let $G$ be a random graph with $V(G)=[n]$ defined by adding each edge $uv$ with $u \in X_i$ and $v \in X_j$ independently with probability $p_{ij}$. Let $\alpha_i = |X_i|/n$ and let
 	$$L = \sqrt{-\sum_{(i,j) \in [k]^2} \alpha_i\alpha_j\log\s{1-p_{ij}}}.$$
 	If $L > \eps$, then 
 \begin{equation}\label{e:tlower}\Prob \s{\had(G) \leq (1 +\eps) L\frac{n}{\sqrt{\log n}}} \geq 1 -\eps.\end{equation}
\end{thm}

The proof of the upper bound, \cref{mainthm}, extends the argument of Axenovich and Steiner~\cite{hadramtwo} on $2$-decompositions $(G_1,G_2)$ of $K_n$. They show that if  $\had(G_1),\had(G_2) \leq O\fs{n}{\sqrt{\log n}}$  then almost all edges of $G_i$ belong to a $(\delta n)$-connected subgraph of $G_i$ for sufficiently small $\delta$.  They then apply a result of Thomason~\cite[Theorem~4.1]{MR1814910} that essentially states that among all $(\delta n)$-connected graphs on $n$ vertices with given density the random graph minimizes the Hadwiger number.

Neither step directly goes through in our application.  We start with a weakening of the first step of the Axenovich-Steiner~\cite{hadramtwo} argument and show that almost all edges of $G_i$ belong to a union of a collection of vertex disjoint  $(\delta n)$-connected subgraphs of $G_i$. We then deduce that if none of the parts has exceedingly large Hadwiger number then for two out of three parts we can after all assume that a single $(\delta n)$-connected subgraph contains almost all of the edges. To lower bound $\had(G_i)$, we need the following extension of the result of Thomason to a setting of \cref{thomasonlower}, where density is not fixed globally, but with respect to a vertex partition. 

\begin{thm}\label{thomasonupper} For all $\eps > 0$ there exists $n_0=n_{0,\ref{thomasonupper}}(\eps)$ such that for all $n \geq n_0$, all graphs $G$ with $\v(G) = n$ and $\kappa(G) \geq \eps n$, and all partitions $(X_1,\ldots,X_k)$ of $V(G)$ the following holds. Let $p_{ij}$  denote the density of edges between $X_i$ and $X_j$ (in particular, $p_{ii}  = \e(G[X_i])/\binom{|X_i|}{2}$) and let $\alpha_i = |X_i|/n$. Let $p_{ij}^* = \min \{p_{ij}, 1 - \eps\}$. If $\alpha_i \geq \eps$ for all $i \in [k]$, then $$\had(G) \geq \frac{n}{\sqrt{\log n}}\s{\sqrt{-\sum_{(i,j) \in [k]^2} \alpha_i\alpha_j\log\s{1-p^*_{ij}}} -\eps}  $$
\end{thm}

The extension is once again technical, but routine. \cref{thomasonupper} reduces the proof of \cref{mainthm} to an optimization problem. This problem, while unappealing at a first glance, allows a short solution.

\subsection*{AI disclosure}
As mentioned above, the proofs of Theorems~\ref{thomasonlower} and~\ref{thomasonupper} are technical extensions of existing proofs in the literature. They were produced by AI following an outline given by the authors, who subsequently verified and polished them. To separate the AI assisted part of the argument from the rest, we relegate these proofs to the appendix.

The remainder of the argument was initially produced and written down without AI assistance. At the final stage of preparation of the paper, AI was used for proofreading, literature search and creating the figures. At this stage AI was also used to extend our results as discussed in Section~\ref{sec:conclusion}.

\paragraph*{Lean formalization.} A Lean formalization of the proof of Theorem~\ref{t:main}, as well as the two additional results discussed in the conclusion was produced by AI. It is available at \url{https://github.com/NeilRahman/HadwigerRamsey}.

\section{Proof of \cref{existencethm}}

In this section, we derive \cref{existencethm} from \cref{thomasonlower}. 

Clearly, it suffices to prove the theorem for even $n$, and so we assume that $n$ is even for the remainder of the proof.
Let $p=5-2\sqrt{6}$.
We assume without loss of generality that $\eps \leq p$ and let $n_0$ be chosen so that $n_0 \geq  2n_{0,\ref{thomasonlower}}(\frac{\eps}{2})$ where $n_{0,\ref{thomasonlower}}(\frac{\eps}{2})$ is as in \cref{thomasonlower}, and so that  \begin{equation}\label{e:logn2}\s{1+\frac{\eps}{2}}\frac{\sqrt{\log n}}{\sqrt{\log (n/2)}} \leq 1 + \eps\end{equation}
for all $n \geq n_0$. 

Construct a random $3$-decomposition $(G_1,G_2,G_3)$ of $K_n$ as described in the introduction. Let $(X,Y)$ be a partition of the vertex set of $K_n$ with $|X|,|Y| = \frac{n}{2}$. Every edge with both ends in $X$ or both ends in $Y$ is assigned to $G_1$ with probability $1-2p$ and to each of $G_2$ and $G_3$ with probability $p$, independently. Every edge with one end in $X$ and the other in $Y$ is assigned to $G_2$ with probability $1/2$ and, otherwise, to $G_3$.  It suffices to show \begin{equation}\label{e:t2}
	\Prob \s{\had(G_i) > (1+\eps)\beta\frac{n}{\sqrt{\log n} }} < \frac{1}{3}	
	\end{equation}
for $i=1,2,3$. 

By the choice of $n_0$ and of the decomposition, the random graph $G_1[X]$ satisfies the conclusion \eqref{e:tlower} of Theorem~\ref{thomasonlower} with $k=1,\alpha_1 = 1,$ $n/2$ replacing $n$, $\eps/2$ replacing $\eps$ and $p_{11}=1-2p$, i.e. $$\Prob\s{\had(G_1[X]) > \s{1+\frac{\eps}{2}}\sqrt{-\log(2p)} \frac{n}{2\sqrt{\log (n/2)}} } \leq \frac{\eps}{2}.$$
As $\sqrt{-\log(2p)}=2\beta$ and \eqref{e:logn2} holds this implies
$$\Prob\s{\had(G_1[X]) > (1+\eps)\beta\frac{n}{\sqrt{\log n} }} \leq \frac{\eps}{2}.$$
By symmetry the same bound holds for $\had(G_1[Y])$. As $G_1$ is a disjoint union of $G_1[X]$ and $G_1[Y]$, we have $\had(G_1) = \max\s{\had(G_1[X]),\had(G_1[Y])}$, and thus
by the union bound \eqref{e:t2} holds for $i=1$.

It remains to verify \eqref{e:t2} for $i=2,3$. The random graph $G_2$ satisfies  the conclusion \eqref{e:tlower} with $k=2,\alpha_1 = \alpha_2=1/2, p_{11}=p_{22}=p, p_{12}=1/2$, and $\eps/2$ in place of $\eps$. Thus 
$$\Prob\s{\had(G_2) > \s{1+\eps/2}\sqrt{-\frac{1}{2}\log(1-p) -\frac{1}{2}\log \frac{1}{2}}\frac{n}{\sqrt{\log n}}} \leq \eps/2.$$
As $\sqrt{-\frac{1}{2}\log(1-p) -\frac{1}{2}\log \frac{1}{2}}=\beta$, and $\eps/2 < 1/3$, this implies that \eqref{e:t2} holds for $i=2$. The case $i=3$ is symmetric.

\section{Proof of \cref{mainthm}.}

Our first major goal  is \cref{connectedsets} which shows that for  every $n$-vertex graph we can delete a small fraction  of  vertices,  and  of possible edges so that the remaining graph is a disjoint union of pieces with connectivity linear in  $n$. (The actual statement is a bit more restrictive.) This lemma is then used in several ways. In Lemma~\ref{twobigsets} \cref{connectedsets} is used to show that if a $3$-decomposition $(G_1,G_2,G_3)$ does not satisfy the conclusion of Theorem~\ref{mainthm} then two of the three parts have almost linear global connectivity, i.e. their restrictions to some pieces containing almost all the vertices have linear connectivity. In the main proof this decomposition of the remaining part is used to reduce it to an optimization problem.

We start with a simple auxiliary lemma, which shows that if two partitions of a finite set do not both contain a part containing most of the vertices, then some two linearly large sets in both partitions individually belong to a single part, but are separated from each other.

\begin{lem}\label[lem]{bigset}
Let $\mathcal{P}_1, \mathcal{P}_2$ be partitions of the same nonempty $n$-element set into $\leq k$ parts. Then for every $\eps >0$ either

\begin{itemize}\item For some $i \in \{1, 2\}$, there exists $P' \in \mathcal{P}_i$ such that $|P'| \geq (1 - \eps)n$, or
\item there exist $A, B$ such that $|A|, |B| \geq \frac\eps{k^2} n$, and for all $i \in \{1,2\}$, there exist $P_{A,i}, P_{B,i} \in \mathcal{P}_i$ such that $A \subseteq P_{A,i}$, $B \subseteq P_{B,i}$, and $P_{A,i} \neq P_{B,i}$ \end{itemize}
\end{lem}

\begin{proof} Suppose the second condition does not hold. If we had two pairs $(P_1, Q_1), (P_2, Q_2) \in \mathcal{P}_1 \times \mathcal{P}_2$ with $|P_i \cap Q_i| \geq \frac\eps{k^2} n$ for $i \in \{1, 2\}$, $P_1 \neq P_2$, $Q_1 \neq Q_2$, we could choose $A = P_1 \cap Q_1$ and $B = P_2 \cap Q_2$ to fit the second condition. Therefore, without loss of generality there exists $P' \in \mc{P}_1$ such that for all $(P, Q) \in \mathcal{P}_1 \times \mathcal{P}_2$ with $|P \cap Q| \geq \frac\eps{k^2} n$, we have $P = P'$. But then, for all $P \in \mathcal{P}_1 \setminus \{P'\}$, we have $$|P| \leq \sum_{Q \in \mathcal{P}_2} |P \cap Q| \leq k  \cdot \frac\eps{k^2}n = \frac\eps k n.$$ Since $|\mathcal{P}_1| \leq k$  it follows that $$n - |P'|  = \sum_{P \in \mathcal{P}_1 \setminus \{P'\}} |P|\leq (k-1) \cdot \frac\eps k n \leq \eps n,$$ and therefore  $|P'| \geq (1-\eps)n$.
\end{proof}

Next we precisely define the properties of the partition that we informally described at the beginning of the section.

Let $G$ be a graph with $|V(G)| = n$. For $\eps,\eta,\delta >0$ and an integer $k \geq 0$
we say that a partition $(V_1, \ldots, V_k, X, Y)$ of $V(G)$  is an \defn{$(\eps,\eta,\delta)$-controlled partition of $G$} if \begin{itemize}
	\item[(a)] $|X| \leq \eps n$,
	\item[(b)] $G[V_i]$ is $(\delta n)$-connected for every $1 \leq i \leq k$,
	\item[(c)] $G[Y]$ is $(3\eps n)$-degenerate,
	\item[(d)] For all $v \in Y$, $v$ has at most $\eta n$ neighbours in $\bigcup_{i=1}^k V_i$,
	\item[(e)] $|V_i| \geq \eps n$ for every $1 \leq i \leq k$,
	\item[(f)] $E(V_i, V_j) = \emptyset$ for all $i \neq j$.
\end{itemize}
Note that if $(V_1, \ldots, V_k, X, Y)$ is a partition of $V(G)$ satisfying (e) then $k \leq \eps^{-1}$. In particular this bound on $k$ holds for every $(\eps,\eta,\delta)$-controlled partition.

\begin{lem}\label[lem]{connectedsets} For all $\eps, \eta > 0$, there exists $\delta = \delta_{\ref{connectedsets}}(\eps,\eta) >0$ such that every graph $G$ admits an $(\eps,\eta,\delta)$-controlled partition.
\end{lem}

\begin{proof}  We assume that $G$ is not $(3\eps n)$-degenerate, as otherwise, the partition $(\emptyset, V(G))$ is $(\eps,\eta,\delta)$-controlled for any choices of $\eta$ and $\delta$. In particular, we assume  $\eps \leq 1/3$. We show that
$$\delta = \min\s{\frac{\eps^3}{(\lfloor\eps^{-1}\rfloor)!},\frac{\eps\eta}{(\lfloor\eps^{-1}\rfloor)!}}$$
satisfies the lemma. 
	
We start by  considering a preliminary partition $(U_1,\ldots, U_{k}, X', Y)$ of $V(G)$ such that \begin{itemize}
	\item [(a')] $|X'| \leq k\eps^{3}n$,
	\item [(b')] $G[U_i]$ is $(\eps^3 n)$-connected for every $1 \leq i \leq k$,
	\item [(c)] $G[Y]$ is $(3\eps n)$-degenerate,
	\item [(e')] $U_i \neq \emptyset$ and $\deg_{G[U_i]}(v) \geq 3\eps n - k\eps^3 n$ for every $1 \leq i \leq k$ and $v \in V(G[U_i])$, and 
	\item[(k)] $k\leq \eps^{-1}$.
 \end{itemize}		
Let us show that it exists. Let $H$ be a maximal subgraph of $G$ with minimum degree $3 \eps n$, then let $Y = V(G) \setminus V(H)$. $G[Y]$ is $(3 \eps n)$-degenerate, as otherwise, there would be some $Y' \subseteq Y$ such that $H \cup G[Y']$ contradicts the choice of $H$.  Therefore the partition $(V(G)-Y, \emptyset, Y)$ satisfies (a'), (c) and (e') with $k=1$.

Choose a partition    $(U_1,\ldots, U_{k'}, X', Y)$ of $V(G)$ satisfying (a'), (c), (e') and (k) with $k'$ maximum. 
We claim that this partition satisfies (b'). 

Suppose not. Then without loss of generality there exists a proper separation $(W,W')$ of $U_{k'}$ with $|W \cap W'| \leq \eps^3 n$. But in such a case the partition  $(U_1,\ldots, U_{k'-1}, W' -W, W-W',  X' \cup (W \cap W'), Y)$ of $V(G)$ contradicts our choice of partition. Indeed, it is easy to see that conditions (a'), (c), (e') are satisfied for the new partition. To verify (k) for the new partition, we need to show that $k' \leq \eps^{-1}-1$. But, the conditions (e') and (k) for the current partition imply that
$$|U_i| \geq  3\eps n - k' \eps^3 n \geq 3\eps n - \eps^2 n  \geq 2 \eps n$$
for every $i$, and so $|V(G)| \geq \sum_{i=1}^{k'}|U_i| \geq 2\eps k' n$, implying $k' \leq \frac{1}{2}\eps^{-1} \leq \eps^{-1}-1,$ where the last inequality follows from $\eps \leq 1/3$, as needed.

This finishes the construction of the preliminary partition. 
We now replace the conditions (a'), (b') and (e') by
 \begin{itemize}
    \item [(a'')] $|X'| \leq \eps^2 n$,
 	\item [(b'')] $G[U_i]$ is $(k! \delta n)$-connected for every $1 \leq i \leq k$,
 	\item [(e'')] $|U_i| \geq 2 \eps n$ for   every $1 \leq i \leq k$.
 \end{itemize}	Note that (e'') holds for the current preliminary partition as shown above, (a') and (k) imply (a''), and  $k! \delta n \leq \eps^3 n$ by (k) and the choice of $\delta$, so (b') implies (b''). 

We now  choose  a partition $(U_1,\ldots, U_{k}, X', Y)$ of $V(G)$ satisfying  (a''),(b''),(c) and (e'') with $|Y|$ minimum. (Note that (k) is implied by (e'') and so also still holds.) Then  \begin{itemize}
    \item[(d')] For all $v \in Y$, $v$ has at most $\eta n$ neighbours in $\bigcup_{i=1}^k U_i$,
\end{itemize} holds. Indeed, if  some $v \in Y$ has at least  $\eta n$ neighbours in $\bigcup_{i=1}^k U_i$, then it has at least $\eps \eta  n$ neighbours in some $U_i$. As $\eps \eta \geq k! \delta$, replacing $Y$ by $Y - \{v\}$ and $U_i$ by $U_i \cup \{v\}$ in this case contradicts minimality of $|Y|$.

Thus we can  choose  a partition $(U_1,\ldots, U_{k}, X', Y)$ of $V(G)$ satisfying  (a''),(b''),(c),(d') and (e'') and we choose it with $k$ minimum. If $k = 0$, it is easy to see that we are done, so suppose that $k \geq 1$. 
Then $G[U_i\cup U_j]$ is not  $((k-1)! \delta n)$-connected  for all $1 \leq i, j \leq k$, $i \neq j$  as otherwise we would replace $U_i$ and $U_j$ by $U_i \cup U_j$ in the partition. It follows that for all such pairs $i,j$ there exists $Z_{ij} = Z_{ji} \subseteq (U_i \cup U_j)$ with $|Z_{ij}| \leq (k-1)!\delta n$ such that $E(U_i-Z_{ij},U_j - Z_{ij})=\emptyset$.  

Finally, define
$$V_i = U_i - \bigcup_{j \neq i}Z_{ij} \qquad \mathrm{and} \qquad  X = X' \cup \s{\bigcup_{\substack{1 \leq i,j \leq k \\ i \neq j}} Z_{ij}}.$$
We have $\sum_{j \neq i}|Z_{ij}| \leq (k-1)(k-1)! \delta n$. Thus $$|V_i| \geq |U_i|-(k-1)(k-1)! \delta n \geq 2\eps n-(k-1)(k-1)! \delta n \geq \eps n.$$  As $G[U_i]$ was $(k!\delta n)$-connected,  the graph $G[V_i]$ is $((k! - (k-1)(k-1)!)\delta n)$-connected and so is $(\delta n)$-connected. Thus (b) and (e) hold. As $$|X| \leq |X'| + k^2(k-1)!\delta n \leq \eps^2n + k^2(k-1)!\delta n \leq \eps n,$$ (a) holds. The condition (c) already holds, (d') implies (d), and (f) holds by the choice of the sets $Z_{ij}$. We conclude that  $(V_1, \ldots, V_k, X, Y)$ is $(\eps,\eta,\delta)$-controlled . 
\end{proof}

Before proceeding to an application of  \cref{connectedsets}, we present two technical consequences of \cref{thomasonupper}. These show that dense graphs and bipartite graphs on $n'$ vertices and with density close to maximum contain complete minors of size $\omega\fs{n'}{\sqrt{\log n'}}$.

\begin{lem}\label[lem]{bipMinor} For every $K > 0$ there exist $\gamma=\gamma_{\ref{bipMinor}}(K), n_0=n_{0,\ref{bipMinor}}(K) >0$ satisfying the following.  
	Let $H$ be a bipartite graph with bipartition $(A,B)$ with $|A|=|B| \geq n_0$ and $|E(H)|\geq (1-\gamma)|A|^2$. Then $$\had(H) \geq K \frac{|A|}{\sqrt{\log{|A|}}}.$$
	
\end{lem}
\begin{proof} We assume $K \geq 1$ without loss of generality. 
	Let $\eps$ be chosen so that $\sqrt{-\log(\eps)} \geq 8 K.$ We show that $\gamma = \frac{\eps}{2} \leq 0.01$ and $n_0 = \max \s{ n_{0,\ref{thomasonupper}}(\eps),100}$ suffice, where $n_{0,\ref{thomasonupper}}(\eps)$ is as in \cref{thomasonupper}.
		 
Let $a = |A|$. Repeatedly delete vertices of $H$ with degree at most $\frac34a$. Suppose we perform $2c$ such deletions for some $c \leq \frac14 a$. We obtain a subgraph of $H$ with $2(a - c)$ vertices and at least $(1-\gamma)a^2-\frac32 ac$ edges. As every $2(a-c)$ vertex subgraph of $H$ has at most $(a-c)^2$ edges we have $(1-\gamma)a^2 -\frac{3}{2}ac \leq (a-c)^2$, i. e.
$$\frac12 ac - \gamma a^2 \leq c^2$$
$$\implies (c-\frac14a)^2 \geq \s{\frac1{16} - \gamma}a^2.$$
Since $c \leq \frac14a$ and $\gamma \leq 0.01$,
$$\frac14a - c \geq \sqrt{\frac1{16}-\gamma} a$$
Given $\gamma \leq 0.01$, we must have $c \leq \frac{1}{10}{a}$. It follows that the process must stop before arriving at the null graph, indeed otherwise the process jumps from at most $\frac{1}{5}a$ deleted vertices to more than $\frac14a$ deleted vertices in two deletions, implying $\frac1{20} a \leq 2$, contradicting $a \geq n_0 \geq 100$.
	
The graph $H'$ obtained after repeated deletions has bipartition $(A',B')$ such that $ |A'|,|B'| \leq a$,  minimum degree at least $\frac{3}{4}a$ and at least $|A'||B'|-\gamma a^2$ edges. In particular, $\kappa(H') \geq \frac{3}{4}a \geq \eps |V(H')|$ and $$\frac{|E(H')|}{|A'||B'|} \geq  1 - \frac{\gamma a^2}{\s{\frac{3}{4}a}^2} \geq 1 -2 \gamma.$$ We apply \cref{thomasonupper} to the partition $(A',B')$ of the vertex set of $H'$ with $$p^*_{12}=\min\s{\frac{|E(H')|}{|A'||B'|}, 1-\eps} =1-\eps$$ to obtain
\begin{align*}\had(H) &\geq \had(H') \geq \frac{|V(H')|}{\sqrt{\log |V(H')|}}\left(\sqrt{-\frac{2|A'||B'|}{(|A'|+|B'|)^2} \log(\eps)}-\eps\right)\\  &\geq \frac{1}{2}\frac{|A|}{2\sqrt{\log {|A|}}}\sqrt{-\frac14\log(\eps)} \geq \frac18\frac{|A|}{\sqrt{\log {|A|}}}\cdot 8K = K \frac{|A|}{\sqrt{\log{|A|}}}.  \end{align*} 
\end{proof}

A stronger version of \cref{bipMinor} easily follows.

\begin{cor}\label[cor]{newBipMinor}
For every $K > 0$ there exist $\gamma=\gamma_{\ref{newBipMinor}}(K), n_0=n_{0,\ref{newBipMinor}}(K) >0$ satisfying the following.  
	Let $H$ be a bipartite graph with bipartition $(A,B)$ with $|B| \geq |A| \geq n_0$ and $|E(H)|\geq (1-\gamma)|A||B|$. Then $$\had(H) \geq K \frac{|A|}{\sqrt{\log{|A|}}}.$$
\end{cor}

\begin{proof}
Let $\gamma = \gamma_{\ref{bipMinor}}(K)$ and $n_0 = n_{0,\ref{bipMinor}}(K)$ be as in \cref{bipMinor}. Let $a = |A|$. Let $B' \subseteq B$ with $|B'| = a$ be chosen uniformly at random. Then $\mathbb{E}|E_H(A,B')| = \frac a{|B|} |E(H)| \geq (1-\gamma) a^2$. Fix a choice $B'$ with $|E_H(A, B')| \geq (1-\gamma)a^2$. By the choice of $\gamma$ and $n_0$, the conclusion of \cref{bipMinor} holds for the bipartite subgraph $H'$ of $H$ with bipartition $(A, B')$ and $E(H') = E_H(A, B')$, so we have
$$\had(H) \geq \had(H') \geq K \frac{a}{\sqrt{\log a}} = K \frac{|A|}{\sqrt{\log{|A|}}},$$
as desired.
\end{proof}

We will use \cref{newBipMinor} for the proof of \cref{mainthm}. However, for the following claim, \cref{bipMinor} suffices.

\begin{lem}\label[lem]{denseMinor} For every $K > 0$ there exist $\gamma=\gamma_{\ref{denseMinor}}(K), n_0=n_{0,\ref{denseMinor}}(K) >0$ satisfying the following.  
Let $H$ be a graph with $|V(H)| \geq n_0$ and $|E(H)|\geq (1-\gamma)\binom{|V(H)|}{2}$. Then $$\had(H) \geq K \frac{|V(H)|}{\sqrt{\log{|V(H)|}}}.$$
\end{lem}

\begin{proof} Let $\gamma= \gamma_{\ref{bipMinor}}(3K)/2$ be as in \cref{bipMinor} and let $n_0 \geq 3n_{0,\ref{bipMinor}}(3K)$ be chosen so that $(1-\gamma)\frac{n_0-1}{n_0} \geq 1 - 2\gamma$.
	Let $a = \lfloor\frac{|V(H)|}{2}\rfloor$. Let a pair of disjoint $A,B \subseteq V(H)$  with $|A|=|B|=a$ be chosen uniformly at random. Then $\bb{E}|E_H(A,B)| \geq \frac{2a^2|E(H)|}{|V(H)|^2} \geq (1-2\gamma)a^2.$ Fix a choice $(A,B)$ with $|E_H(A,B)| \geq (1-2\gamma)a^2$. By the choice of $\gamma$ and $n_0$, the conclusion of \cref{bipMinor} with $3K$ in place of $K$ holds for the bipartite subgraph $H'$ of $H$ with bipartition $(A,B)$ and $E(H')=E_H(A,B)$, and so we have $$\had(H) \geq \had(H') \geq 3K \frac{a}{\sqrt{\log{a}}} \geq K \frac{|V(H)|}{\sqrt{\log{|V(H)|}}},$$
	as desired.
\end{proof}


The following is our first application of \cref{connectedsets} which also uses \cref{bigset}. It shows that if a $3$-decomposition $(G_1,G_2,G_3)$ of a complete graph $K_n$ is not too far from being extremal for Hadwiger-Ramsey numbers then two out of three graphs in the decompositions are either very sparse or have connectivity linear in $n$ after possibly deleting a small exceptional set of vertices.                    

\begin{lem}\label[lem]{twobigsets} For all $\eps ,C > 0$ there exist $n_0 = n_{0,\ref{twobigsets}}(\eps,C), \delta = \delta_{\ref{twobigsets}}(\eps,C) > 0$ such that the following holds.
	
Let $(G_1,G_2,G_3)$ be a  $3$-decomposition of a complete graph $G =K_n$ for $n \geq n_0$.  Then either 
\begin{itemize}
\item There exist two distinct $i_1,i_2 \in [3]$ such that for each $j \in [2]$ \begin{itemize} \item either  $|E(G_{i_j})| \leq 3\eps n^2$, \item or there exists $P_j \subseteq V(G)$ such that $|P_j| \geq (1-\eps)n$ and $G_{i_j}[P_j]$ is $(\delta n)$-connected,\end{itemize} or
\item  $\had(G_i) \geq C\frac n{\sqrt{\log n}}$ for some $i \in [3]$.
\end{itemize}
\end{lem}

\begin{proof} Let $0<\eps'<\frac1{10}$ be chosen so that 
$
(1-\eps')(1-2\eps') > (1-\eps)$ and $\eps' \leq \frac{2}{3}\eps$.
 Let $K=2C(\eps')^{-3}$. Let $\gamma=\gamma_{\ref{bipMinor}}(K)$ be as in \cref{bipMinor} and let $n'_0=\max(n_{0,\ref{bipMinor}}(K), 2)$. Let $\eta =\frac{1}{4}(\eps')^{3}\gamma$.
By \cref{connectedsets} there exists $\delta > 0$ such that every graph admits an $(\eps',\eta,\delta)$-controlled partition. Finally, let $n_0 = 2 (\eps')^{-3}n'_0$. We show that these choices of $\delta$ and $n_0$ satisfy the lemma.

Assume that the first outcome of the lemma does not hold, i.e.\ we assume without loss of generality that $|E(G_{i})| \geq 3\eps n^2$  and there does not exist $P_i \subseteq V(G)$ such that $|P_i| \geq (1-\eps)n$ and $G_{i}[P_i]$ is $(\delta n)$-connected for $i=1,2$. 
Let 
$(V_{1,i}, \ldots, V_{k_i,i}, X_i, Y_i)$ be  an $(\eps',\eta,\delta)$-controlled partition of $G_i$. As $G_i[V_{j,i}]$ is $(\delta n)$-connected, we have $|V_{j,i}| \leq (1-\eps)n$ for  $i=1,2$ and $j \in [k_i]$. Moreover, $$3\eps n^2 \leq |E(G_i)| \leq |E(G_i[Y_i])| + |V(G)-Y_i| n \leq 3\eps'n^2 + |V(G)-Y_i| n,$$ implying $|V(G)-Y_i| \geq 3(\eps -\eps')n \geq \eps n$.

Apply \cref{bigset} with $\eps'$ in place of $\eps$ and $k =\max(k_1,k_2)+1$ to the pair of  partitions
$$\mathcal{P}_1 = (V_{1,1} \setminus X_2, \ldots, V_{k_1,1} \setminus X_2, Y_1 \setminus X_2),$$
$$\mathcal{P}_2 = (V_{1,2} \setminus X_1, \ldots, V_{k_2,2} \setminus X_1, Y_2 \setminus X_1)$$
of $V(G) \setminus X_1 \setminus X_2$. 
As we have shown above, each part of $\mc{P}_i$ has size at most $$(1-\eps)n < (1-\eps')(1-2\eps')n \leq (1-\eps')|V(G) \setminus X_1 \setminus X_2|.$$ Thus the first outcome of \cref{bigset} cannot hold.
Further note that $k \leq (\eps')^{-1}+1$.
 So there exist $A, B \subseteq V(G)$ such that
$$|A|, |B| \geq \frac{\eps'}{(k+1)^2}(1-2\eps') n \geq \frac{\eps'(1-2\eps')}{((\eps')^{-1}+2)^2} n\geq \frac{1}{2}(\eps')^3n$$
and for all $i \in \{1,2\}$, there exist $P_{A,i}, P_{B,i} \in \mathcal{P}_i$ such that $A \subseteq P_{A,i}$, $B \subseteq P_{B,i}$, and $P_{A,i} \neq P_{B,i}$. Without loss of generality we assume $|A|=|B|$. Let $a = |A|$. 
By definition of the $(\eps',\eta,\delta)$-controlled partition it follows that  either every vertex in $P_{A,i}$ has at most $\eta n$ neighbours in $P_{B,i}$ or every vertex of 
 $P_{B,i}$ has at most $\eta n$ neighbors in $P_{A,i}$. Thus there are at most $\eta n a$ edges between $A$ and $B$ in each of $G_1$ and $G_2$. 
 
Consider the bipartite subgraph $H$ of $G_3$ with $V(H)=A \cup B$ and $E(H)=E_{G_3}(A,B)$. We have $$|E(H)| \geq a^2 - 2\eta n a \geq (1- 4(\eps')^{-3}\eta) a^2 = (1-\gamma)a^2$$ and $a \geq \frac{1}{2}(\eps')^3 n \geq n'_0$. By \cref{bipMinor},
$$\had(H) \geq K\frac{ a}{\sqrt{\log a}} \geq  K\frac {\frac{1}{2}(\eps')^3n}{\sqrt{\log n}} \geq C\frac{ n}{\sqrt{\log n}},$$
and so the second outcome of the lemma holds with $i=3$.
\end{proof}

The next two lemmas solve the optimization problem which the proof eventually reduces to.

\begin{lem}\label[lem]{fopt} Let $$f(x) =   x^2 \log\s{2 - e^{-\beta^2/x^2}} \quad \mathrm{for} \quad 0<x\leq1 \qquad \mathrm{and} \qquad f(0)=0. $$
Let $\alpha_1, \ldots, \alpha_k \geq  0$ with $\sum_{i=1}^k \alpha_i = 1$. Then
$$ \sum_{i=1}^{k} f(\alpha_i) \leq \log 2 - \beta^2.$$	
\end{lem}	

\begin{proof} Let $F = \log 2 - \beta^2$. We have $2f(1/2)=F$.
	The function $f(x)/x$ increases on the interval $(0,1/2]$. If $\alpha_i < 1/2$ for every $i$, then $$\sum_{i=1}^{k} f(\alpha_i) \leq \sum_{i=1}^k \alpha_i \frac{f(1/2)}{1/2} = F$$ and the lemma holds.
	Otherwise, assume without loss of generality  that  $\alpha_1 \geq 1/2$. If $\alpha_1 = 1$, $\alpha_i = 0$ for all $i > 1$, so the sum is $f(1)=0.284\ldots < 0.293\ldots=F$. If instead $\alpha_1 < 1$, $$\sum_{i=2}^{k} f(\alpha_i) \leq \sum_{i=2}^k \alpha_i \frac{f(1-\alpha_1)}{1-\alpha_1} = f(1-\alpha_1),$$
	and it remains to show that $f(\alpha_1)+f(1-\alpha_1) \leq F$, i.e.\ to prove the lemma for $k=2$. Let $h(x)=f(x)+f(1-x).$ The function $h'(x)$  has three roots on $[0,1]$ at approximately $x_0=0.165..., 1/2,$ and $1-x_0$, and $h(x)$ decreases on $[0,x_0]$ and increases on $[x_0,1/2]$. Thus $h(x)$ is maximized at $x=0$ or $x=1/2$. We have  $h(1/2)=F=0.293\ldots$ and $h(1)=0.284\ldots < F$.      
\end{proof}

\begin{lem}\label[lem]{opt} Let $\alpha_1, \ldots, \alpha_\ell > 0$ with $\sum_{i=1}^\ell \alpha_i = 1$. Let $ 0 \leq p_1, \ldots, p_\ell < 1$. Let $ 0 \leq q_{ij}^{B}, q_{ij}^G <1$ for all $(i,j) \in [\ell]^2$. If \begin{align*} &p_i + q_{ii}^B + q_{ii}^G = 1 \qquad \mathrm{for\; all} \quad i \in [\ell], \qquad \mathrm{and} \\ &q_{ij}^B + q_{ij}^G = 1  \qquad \mathrm{for\; all} \quad i, j \in [\ell] \quad \mathrm{with} \; i \neq j,\end{align*} then at least one of the following holds:
\begin{itemize}
\item $\alpha_i \sqrt{-\log\s{1-p_i}} \geq \beta$ for some $i \in [\ell]$, or
\item $\sqrt{-\sum_{(i,j) \in [\ell]^2} \alpha_i\alpha_j\log\s{1-q^B_{ij}}} \geq \beta$, or
\item $\sqrt{-\sum_{(i,j)\in [\ell]^2} \alpha_i\alpha_j\log\s{1-q^G_{ij}}} \geq \beta$.
\end{itemize}
\end{lem}

\begin{proof} Suppose not, i.e.\ after squaring we have
\begin{equation}\label{e:first} -\alpha^2_i \log\s{1-p_i} < \beta^2 \end{equation}  for all $i \in [\ell]$,
$$-\sum_{(i,j)\in [\ell]^2} \alpha_i\alpha_j\log\s{1-q^B_{ij}}< \beta^2 \quad \mathrm{and} \quad
-\sum_{i,j \in [\ell]} \alpha_i\alpha_j\log\s{1-q^G_{ij}} < \beta^2. $$
Averaging the last two inequalities gives
$$- \sum_{i,j \in [\ell]} \alpha_i \alpha_j \frac{\log\s{1 - q_{ij}^B} + \log\s{1-q_{ij}^G}}2 < \beta^2.$$
By convexity of  $-\log\s{1-x}$ for $x < 1$, we have
$$\sum_{(i,j) \in [\ell]^2} - \alpha_i \alpha_j \log\s{1 - \frac{q_{ij}^B + q_{ij}^G}2} \leq \sum_{(i,j) \in [\ell]^2} - \alpha_i \alpha_j \frac{\log\s{1 - q_{ij}^B} + \log\s{1-q_{ij}^G}}2.$$
Thus $$\sum_{i,j \in [\ell]} - \alpha_i \alpha_j \log\s{1 - \frac{q_{ij}^B + q_{ij}^G}2} < \beta^2.$$ As $q_{ii}^B + q_{ii}^G = 1 - p_i$ and $q_{ij}^B + q_{ij}^G = 1$ for $ i \neq j$ we rewrite the above $$\log 2 \cdot \sum_{i \neq j} \alpha_i \alpha_j - \sum_{i \in [\ell]} \alpha_i^2 \log\s{1 - \frac{1 - p_i}2} < \beta^2.$$ Meanwhile, by \eqref{e:first} we have
$$ (1 - p_i)/2 \geq e^{-\beta^2/\alpha_i^2}/2.$$
Thus
\begin{align*} \beta^2 &> \log2 \cdot \sum_{i \neq j} \alpha_i \alpha_j - \sum_{i \in [\ell]} \alpha_i^2 \log\s{1 - e^{-\beta^2/\alpha_i^2}/2} \\
&=\log 2\cdot\s{1 -\sum_{i \in [\ell]} \alpha_i^2} - \sum_{i \in [\ell]} \alpha_i^2 \log\s{1 - e^{-\beta^2/\alpha_i^2}/2} \\ &
=\log 2 - \sum_{i \in [\ell]} \alpha_i^2 \s{\log2 +\log\s{1 - e^{-\beta^2/\alpha_i^2}/2}} \\ &=
\log 2 - \sum_{i \in [\ell]} \alpha_i^2 \log\s{2 - e^{-\beta^2/\alpha_i^2}}.\end{align*}
Let
$f(x) = x^2 \log\s{2 - e^{-\beta^2/x^2}}$ be as in \cref{fopt}. Then the last inequality can be rewritten as $\sum_{i=1}^\ell f(\alpha_i) > \log 2 - \beta^2$, contradicting \cref{fopt}.
\end{proof}

As all the functions in \cref{opt} are continuous when the variables are bounded away from one, compactness immediately implies the following approximate version of \cref{opt}. 

\begin{cor}\label[cor]{optapprox} For every $\zeta > 0$ there exists $\delta=\delta_{\ref{optapprox}}(\zeta) > 0$ satisfying the following. Let $\ell \leq 1/\zeta $ be a positive integer, let $\alpha_1, \ldots, \alpha_\ell \geq 0$ with $\sum_{i=1}^\ell \alpha_i = 1$. Let $ 0 \leq p_1, \ldots, p_\ell \leq 1 - \zeta $. Let $ 0 \leq q_{ij}^{B}, q_{ij}^G \leq 1 -\zeta $  with $q^{B}_{ij} = q^{B}_{ji},q^{G}_{ij} = q^{G}_{ji}$ for all $(i,j) \in [\ell]^2$. Then the following conditions cannot all simultaneously hold:
	\begin{itemize}
	 \item[(a)] $p_i + q_{ii}^B + q_{ii}^G  \geq  1 - \delta$ for all $i \in [\ell]$,
	 \item[(b)]  $q_{ij}^B + q_{ij}^G  \geq 1 - \delta$  for all $i, j \in [\ell]$ with $i \neq j$, 
		\item[(c)] $\alpha_i \sqrt{-\log\s{1-p_i}} \leq \beta - \zeta $ for all $i \in [\ell]$, 
		\item[(d)] $\sqrt{-\sum_{(i,j) \in [\ell]^2} \alpha_i\alpha_j\log\s{1-q^B_{ij}}} \leq \beta -\zeta $,
		\item [(e)] $\sqrt{-\sum_{(i,j)\in [\ell]^2} \alpha_i\alpha_j\log\s{1-q^G_{ij}}} \leq \beta - \zeta $.
	\end{itemize}
\end{cor}

\begin{proof}[Proof of \cref{mainthm}]
We first introduce the parameters that we will need. We assume $\eps \leq 1$ without loss of generality. Choose $\lambda < \frac{\eps \beta^2}{420 \log 2}$ then let $K = \frac{2\beta}{\lambda}$. Let $\gamma = \min\s{\gamma_{\ref{newBipMinor}}(K), \gamma_{\ref{denseMinor}}(K)}$, let $\zeta = \min\s{\gamma, \frac\lambda2, \frac\beta{100}}$, and let $\xi = \delta_{\ref{optapprox}}(\zeta)$. Choose $\eta < \frac{\lambda \xi}2$ and $\eps_0 = \min\s{\lambda \xi, \lambda \eps, \zeta \beta^2}/100$. Let $\delta_0 = \delta_{\ref{twobigsets}}(\eps_0, \beta)$, let $\delta = \delta_{\ref{connectedsets}}(\lambda, \eta)$, and let $\theta = \min\s{\delta_0/2, \delta, \zeta}$. Finally, it will be clear that we can make $n_0$ sufficiently large based on the other parameters. 
For any $n \geq n_0$, let $(G_1, G_2,G_3)$ be a 3-decomposition of $G =K_n$ and suppose for a contradiction that
\begin{equation}\label{e:lower2} \had(G_i) < (1- \eps) \beta \frac{n}{\sqrt{\log n}}, \end{equation}
for all $i \in [3]$. By \cref{twobigsets}  applied with parameters $\eps_0$  and $C=\beta$, we assume without loss of generality that for $j \in \{1,2\}$  \begin{itemize} \item either $|E(G_j)| \leq 3\eps_0n^2$, or \item 
there exists $P_j \subseteq V(G)$ such that $|P_j|\geq (1-\eps_0)n$ and $G_j[P_j]$ is $(\delta_0 n)$-connected. 
\end{itemize} If either $P_j$ does not exist for $j \in [2]$, then let $P_j = V(G)$. Let $(V_1,V_2,\ldots,V_{k},X,Y)$ be a $(\lambda,\eta,\delta)$-controlled partition of $G_3$. If $|Y| \geq 20\lambda n$, let $V_{k+1}=Y,V_{k+2}=X$ and $\ell=k+2$, otherwise, let $V_{k+1}=X \cup Y$ and $\ell=k+1$. Note that
$$\ell \leq k+2 \leq 2\lambda^{-1} \leq \zeta^{-1}.$$

We  define parameters $\alpha_i, p_i$ for $i \in [\ell]$ and  $q_{ij}^B, q_{ij}^G$ for $(i, j) \in [\ell]^2$. $\alpha_i$ corresponds to the relative size of $V_i$ and $p_i$ corresponds to the edge density of $V_i$ in $G_3$. Within $P_1$, $q_{ij}^B$ corresponds to the $G_1$ edge densities between $V_i$ and $V_j$, and $q_{ij}^G$ is defined analogously with $P_2$ and $G_2$. We use \eqref{e:lower2} together with \cref{thomasonupper} to show that they contradict \cref{optapprox}. 

Let $\alpha_i = \frac{|V_i|}{n}$ for $i \in [\ell]$, \begin{align*}&p_i = \frac{2|E(G_3[V_i])|}{|V_i|(|V_i|-1)}, \qquad q^{B}_{ii} = \frac{|E(G_1[V_i \cap P_1])|}{\binom{|V_i \cap P_1|}2}, \qquad q^{G}_{ii} = \frac{|E(G_2[V_i \cap P_2])|}{{\binom{|V_i \cap P_2|}2}} &\quad \mathrm{for} \; i \in [\ell-1], \; \mathrm{and} \\
 &q^{B}_{ij} = \frac{|E_{G_1}(V_i \cap P_1,V_j \cap P_1)|}{|V_i \cap P_1||V_j \cap P_1|}, \qquad  q^{G}_{ij} = \frac{|E_{G_2}(V_i \cap P_2,V_j \cap P_2)|}{|V_i \cap P_2||V_j \cap P_2|} &\quad \mathrm{for}\;    i,j \in [\ell-1], i \neq j.\end{align*}
Finally,  let $p_{\ell}=0$ and let $q^{B}_{i\ell} = q^{B}_{\ell i} = q^{G}_{i\ell} = q^{G}_{\ell i} = 1/2$ for $i \in [\ell]$. 

If $p_{i} \geq 1-\zeta$ for some $i$ then  by the choices of $\zeta$, $n_0$ and $\lambda$, \cref{denseMinor} is applicable to $G_3[V_i]$ and we have $$\had(G_3) \geq \had(G_3[V_i]) \geq K\frac{|V_i|}{\sqrt{\log |V_i|}} \geq K\lambda \frac{n}{\sqrt{\log n}} \geq \beta  \frac{n}{\sqrt{\log n}},$$
contradicting \eqref{e:lower2}. Thus $p_{i} \leq 1-\zeta$ and similarly $q^{B}_{ij},q^{G}_{ij} \leq 1-\zeta$ for $i,j \in [\ell]$, where for $i \neq j$, $i,j \in [\ell-1]$ we apply \cref{newBipMinor} instead of \cref{denseMinor}. In particular, for $i \in [\ell-1]$ and $j \in [2]$, $|V_i \cap P_j| \geq (\lambda - \eps_0)n \geq \frac\lambda2 n$ by the choice of $\eps_0$, so we can repeat a similar argument to above.

It remains to show that conditions (a)--(e) of \cref{optapprox} are satisfied to obtain a contradiction. (a)--(c) clearly hold when $i = \ell$ or $j = \ell$, so for these conditions, we assume $i \neq \ell$ and $j \neq \ell$.
For (a), $p_i$ is the edge density of $G_3[V_i]$. If we write $q_{ii}^{*B}$ and $q_{ii}^{*G}$ for the  edge densities of $G_1[V_i]$, $G_2[V_i]$, note that $p_i + q_{ii}^{*B} + q_{ii}^{*G} = 1$. Since $|P_j| \geq (1-\eps_0)n$ for $j \in [2]$,
\begin{align*}
q^{B}_{ii} & = \frac{|E(G_1[V_i \cap P_1])|}{\binom{|V_i \cap P_1|}2} \\
& \geq \frac{|E(G_1[V_i])|}{\binom{|V_i|}2} - \frac{ \eps_0n(|V_i|-1)}{\binom{|V_i|}2} \\
&\geq \frac{|E(G_1[V_i])|}{\binom{|V_i|}2} - \frac{2\eps_0}{\lambda}.
\end{align*}
We also have the analogue for $q_{ii}^G$. By the choice of $\eps_0$,
$$p_i + q_{ii}^B + q_{ii}^G \geq 1 - \xi.$$
For (b), we proceed similarly: 
\begin{align*}
q^{B}_{ij} &= \frac{|E_{G_1}(V_i \cap P_1,V_j \cap P_1)|}{|V_i \cap P_1||V_j \cap P_1|} \\
& \geq \frac{|E_{G_1}(V_i,V_j)|}{|V_i||V_j|} - \frac{\eps_0 n (|V_i| + |V_j|)}{|V_i||V_j|} \\
& \geq \frac{|E_{G_1}(V_i,V_j)|}{|V_i||V_j|} - \frac{2 \eps_0}\lambda.
\end{align*}
It remains to verify that $|E_{G_3}(V_i,V_j)|$ is sufficiently small for $i, j < \ell$. Indeed, either $|E_{G_3}(V_i,V_j)| = 0$ if $Y \not\in \{V_i, V_j\}$, or
$$|E_{G_3}(V_i,Y)| \leq \eta n |Y| \implies \frac{|E_{G_3}(V_i,Y)|}{|V_i||Y|} \leq \frac{\eta}{\lambda}.$$
By the choices of $\eta$ and $\eps_0$,
$$q_{ij}^B + q_{ij}^G \geq 1 - \xi.$$
For (c), 
if $V_{\ell-1} = Y$, then $G_3[Y]$ is $(3 \lambda n)$-degenerate, so $|E(G_3[Y])| \leq 3 \lambda n |Y|$. Therefore, $$p_{\ell-1} \leq \frac{6 \lambda n}{|Y| - 1} \leq \frac{6 \lambda n}{20 \lambda n - 1} = \frac3{10} + o(1).$$
Since $\sqrt{-\log\s{1-\frac3{10}}} < 0.99\beta$,
$$\alpha_{\ell-1} \sqrt{-\log(1-p_{\ell-1})} \leq \beta - \zeta$$
by the choice of $\zeta$, for sufficiently large $n$. 
If instead 
$V_i \neq Y$, then $G_3[V_i]$ is $(\delta n)$-connected. By the choice of $\theta$, we can apply \cref{thomasonupper} to $G_3[V_i]$ with parameter $\theta$ and the trivial partition $(V_i)$; in particular, $p^*_{11} = p_i$. We get
$$\had(G_3) \geq \had(G_3[V_i]) \geq \frac{|V_i|}{\sqrt{\log |V_i|}} \s{\sqrt{-\log\s{1-p_i}} - \theta}.$$
By our main assumption \eqref{e:lower2} and the choices of $\zeta$ and $\theta$,
\begin{align*}
& \qquad\:\: \frac{|V_i|}{\sqrt{\log |V_i|}} \s{\sqrt{-\log\s{1-p_i}} - \theta} < (1- \eps) \beta \frac{n}{\sqrt{\log n}} \\
&\implies \alpha_i \sqrt{-\log\s{1-p_i}} < (1 - \eps) \beta \sqrt\frac{\log|V_i|}{\log n} + \alpha_i \theta \\
&\implies \alpha_i \sqrt{-\log\s{1-p_i}} < (1 - \eps) \beta + \alpha_i \theta \\
&\implies \alpha_i \sqrt{-\log\s{1-p_i}} < (1 - \frac\eps2) \beta \\
&\implies \alpha_i \sqrt{-\log\s{1-p_i}} < \beta - \zeta.
\end{align*}
For (d) and (e), the argument is symmetric, so it only remains to prove (d). We can rewrite what we want to show as
$$\s{-\sum_{i, j < \ell} \alpha_i \alpha_j\log\s{1-q_{ij}^B}} + (2 \alpha_\ell - \alpha_\ell^2) \log2 \leq (\beta - \zeta)^2.$$
Since $|X| \leq \lambda n$ and $|Y| \leq 20 \lambda n$ if $Y \subseteq V_\ell$, $2\alpha_\ell \leq 42 \lambda$. Recall that $42 \lambda \log 2 \leq \frac{\eps \beta^2}{10}$. Thus, it suffices to show that
\begin{equation}\label{eq:cased}
-\sum_{i, j < \ell} \alpha_i \alpha_j \log\s{1-q_{ij}^B} < \s{1-\frac\eps2}\beta^2
\end{equation}
We first consider the case where $|E(G_1)| \leq 3 \eps_0 n^2$ and $P_1 = V(G)$. In particular,
$$\sum_{i, j < \ell} \alpha_i \alpha_j q_{ij}^B \leq \frac{2|E(G_1)|}{n^2} + \frac1n \leq 6 \eps_0 + \frac 1n \leq 7 \eps_0$$
for sufficiently large $n$. Note that, as $\frac{-\log(1-x)}x$ is increasing on $(0, 1)$ and $q_{ij}^B \leq 1 - \zeta$, we have $-\log(1-q_{ij}^B) \leq \frac{-\log(\zeta)}{1-\zeta} q_{ij}^B$ for all $i, j$. Therefore,
\begin{align*}
- \sum_{i, j < \ell} \alpha_i \alpha_j \log(1-q_{ij}^B) &\leq -\frac{\log(\zeta)}{1-\zeta} \sum_{i, j < \ell} \alpha_i \alpha_j q_{ij}^B \\
&\leq -\frac{7 \eps_0 \log(\zeta)}{1-\zeta} \\
&\leq -\frac{7\zeta\log(\zeta)}{100(1-\zeta)}\beta^2 < (1 - \frac \eps2) \beta^2
\end{align*}
by the choice of $\eps_0$, establishing \eqref{eq:cased} in this case. Now suppose that $G_1[P_1]$ is $(\delta_0 n)$-connected. If $|V_\ell \cap P_1| \geq (\delta_0/2) n$, let $P' = P_1$ and $\ell' = \ell$. Otherwise, let $P' = P_1 \setminus V_\ell$ and $\ell' = \ell-1$. Note that $G_1[P']$ is still at least $((\delta_0/2) n)$-connected. We apply \cref{thomasonupper} to the partition $(V_1 \cap P', \ldots, V_{\ell'} \cap P')$ of $P'$, with parameter $\theta$, to deduce
$$\had(G_1[P']) \geq \frac{|P'|}{\sqrt{\log |P'|}} \s{\sqrt{-\sum_{i, j < \ell} \alpha_i^* \alpha_j^* \log(1-q^B_{ij})} - \theta},$$
where for $i \in [\ell']$,
\begin{align*}
\alpha^*_i = \frac{|V_i \cap P'|}{|P'|} \geq \min\s{\lambda-\eps_0, \delta_0/2} & \geq \min\s{\lambda/2, \delta_0/2} \\
& \geq \min\s{\zeta, \theta} \geq \theta
\end{align*}
and $q_{ij}^B \leq 1 - \zeta \leq 1 - \theta$. We have that $|V_i \cap P'| \geq (1-\frac{\eps_0}\lambda) |V_i|$ for $i < \ell$, so
$$\had(G_1) \geq \frac{|P'|}{\sqrt{\log |P'|}} \s{\s{1 - \frac{\eps_0}{\lambda}} \frac{n}{|P'|} \sqrt{-\sum_{i, j < \ell} \alpha_i \alpha_j \log\s{1-q^B_{ij}}} - \theta}$$
$$\implies \had(G_1) \geq \frac{n}{\sqrt{\log n}} \s{\s{1 - \frac{\eps_0}{\lambda}} \sqrt{-\sum_{i, j < \ell} \alpha_i \alpha_j \log\s{1-q^B_{ij}}} - \theta}.$$
By our main assumption \eqref{e:lower2},
$$\s{1-\frac{\eps_0}\lambda} \sqrt{-\sum_{i, j < \ell} \alpha_i \alpha_j \log(1-q_{ij}^B)} - \theta < (1-\eps) \beta$$
$$\implies \sqrt{-\sum_{i, j < \ell} \alpha_i \alpha_j \log(1-q_{ij}^B)} < \s{(1-\eps)\beta + \theta}\s{1 - \frac{\eps_0}\lambda}^{-1}$$
$$\leq \frac{1-9\eps/10}{1-\eps/10} \beta$$
$$\leq \sqrt{\s{1-\frac\eps2}}\beta,$$
proving \eqref{eq:cased}, and thus completing the proof.
\end{proof}

\section{Conclusion}\label{sec:conclusion}

In this paper we determined the asymptotic of $3$-color Hadwiger-Ramsey numbers $R_{\had}(k; 3)$. We used Theorems~\ref{thomasonlower}~and~\ref{thomasonupper}, extensions of Bollob\'as, Catlin and Erd\H{o}s upper and Thomason lower bounds on Hadwiger numbers, together with a structural Lemma~\ref{connectedsets} to reduce the problem to an exponential optimization  problem, which we solved in Lemma~\ref{opt}.

We believed that a similar approach could extend to more colors, but as we expected the complexity to grow sharply with the number of colors, we abstained from pursuing this further.

During the final stages of paper revision, as due diligence, we asked AI to extend our approach to four colors. AI succeeded. In particular, it reduced the problem of determining the asymptotics of Hadwiger-Ramsey numbers $R_{\had}(k; \ell)$ for fixed $\ell$ and $k \to \infty$ to an optimization problem as follows. 

Let $\ell$ be the number of colors. Let $k$ be a positive integer and let  $\mc P_1, \ldots \mc P_\ell$ be partitions of $[k]$. Let the parameters $\alpha_i, p^{c}_{ij} \in \bb{R}$ for $i,j \in [k]$ and $c \in [\ell]$ of our optimization problem satisfy the following:   
\begin{itemize}
\item $\alpha_i > 0$ for all $i \in [k]$
\item $\sum_{i=1}^k \alpha_i = 1$
\item $p^c_{ij} = p^c_{ji} \in [0, 1)$ for all $c \in [\ell]$ and all $i, j \in [k]$
\item $p^c_{ij} = 0$ if $i$ and $j$ are in different parts of $\mc P_c$
\item $\sum_{c=1}^\ell p^c_{ij}=1$ for all $i, j \in [k]$
\end{itemize}
Combinatorially, these restrictions correspond to a partition $(X_1,\ldots,X_k)$ of $[n]$ and an $\ell$-decomposition $(G_1,\ldots,G_{\ell})$ of $K_n$ such that $|X_i| = \alpha_i n$, the density of edges of $G_c$ between $X_i$ and $X_j$ is $p^{c}_{ij}$ and the restriction of $G_c$ to $\cup_{i \in P}X_i$ for every $P \in \mc{P}_c$ is sufficiently connected for application of  Theorem~\ref{thomasonupper}.
Define
$$
 \beta_\ell=\inf\;\max_{\substack{c\in[\ell]\\ C\in\mc P_c}}
       \sqrt{-\sum_{i,j\in C}\alpha_i\alpha_j\log(1-p^c_{ij})},
$$
where the infimum is over all instances that fit the above conditions.  

The following result was proved by AI and formalized in Lean.

\begin{prop}\label{conj:optimization}
For any $\ell \geq 2$,
$$R_{\had}(k;\ell) = \s{\beta_\ell^{-1} + o(1)}k \sqrt{\ln k}.$$
\end{prop}

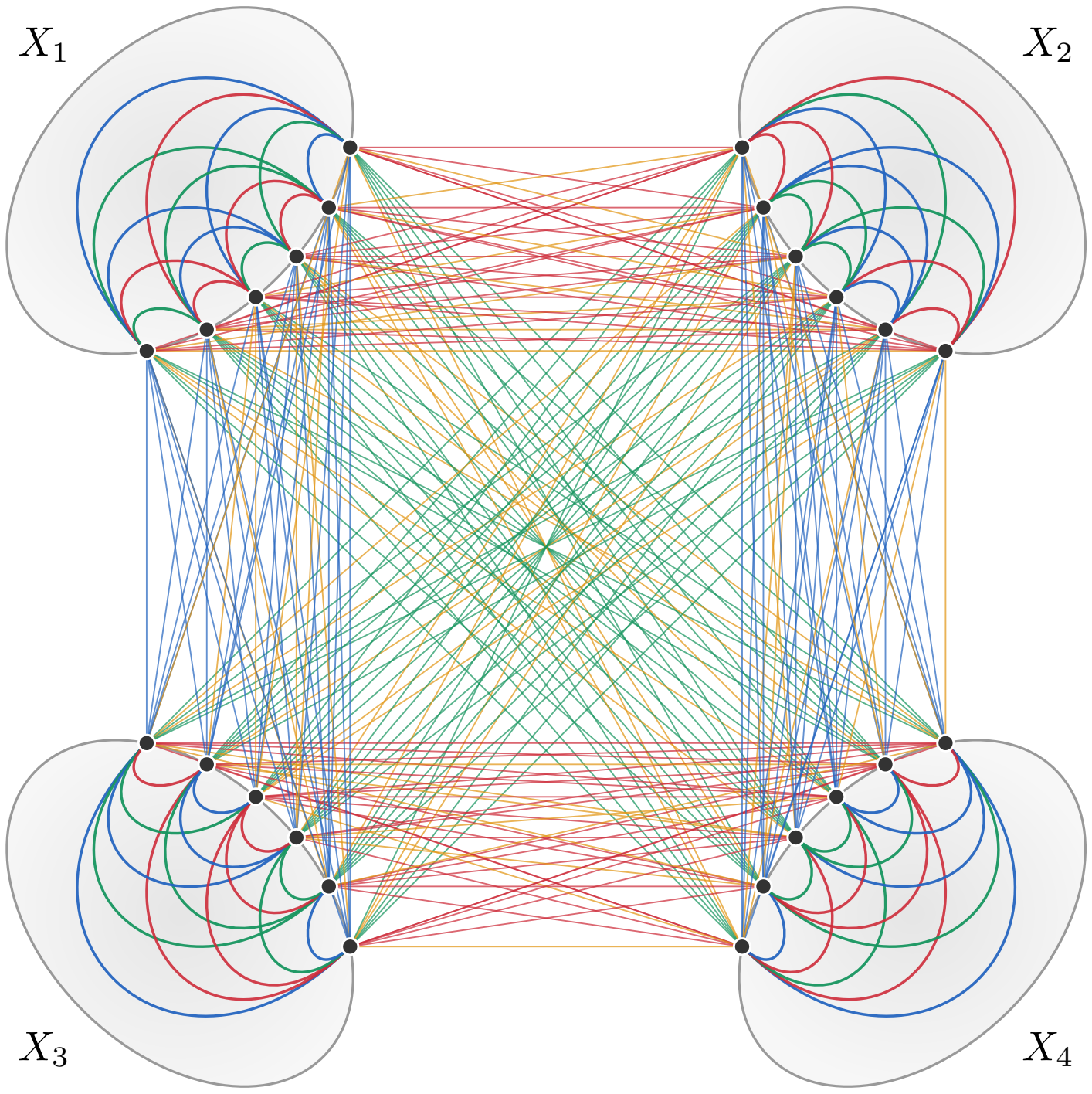
\begin{figure}[t]
	\centering
	\input{four-color.tikz}
	\caption{The optimal construction for the four-color Hadwiger-Ramsey numbers.}
	\label{f:fourcolor}
\end{figure}

For four colors, we initially produced the following construction without AI assistance, before considering the generalized optimization program suggested above. Partition the vertex set of $K_n$ into four equal-size parts $X_1,X_2,X_3,X_4$,  the $4$-decomposition $(G_1,G_2,G_3,G_4)$ is then defined as follows. Every edge is independently at random assigned to one of $G_i$ so that the following is satisfied: 
\begin{itemize}
	\item every edge with both ends in $X_i$ is assigned to each of $G_1,G_2$ and $G_3$ with probability $1/3$;
	\item every edge with one end in $X_1$ and the other in $X_2$, or one end in $X_3$ and the other in $X_4$ is assigned to $G_1$ with probability $p$ and otherwise to $G_4$;
		\item every edge with one end in $X_1$ and the other in $X_3$, or one end in $X_2$ and the other in $X_4$ is assigned to $G_2$ with probability $p$ and otherwise to $G_4$;
	\item every edge with one end in $X_1$ and the other in $X_4$, or one end in $X_2$ and the other in $X_3$ is assigned to $G_3$ with probability $p$ and otherwise to $G_4$.	
\end{itemize} 
See Figure~\ref{f:fourcolor}. Thus the graphs $G_1,G_2$ and $G_3$ are structured similarly, each consisting of two random components of size $n/2$, while the graph $G_4$ has one global component.
Optimization, which we initially performed numerically, suggested $p \approx 0.7445$ minimizes the maximum Hadwiger number of the parts of the decomposition.

AI confirmed optimality of this construction, proving the following Lean formalized proposition. 

\begin{prop}\label{conj:fourcol} $$R_{\had}(k;4) = \s{\beta_4^{-1} + o(1)}k \sqrt{\ln k},$$
where	
$\beta_4 = \sqrt{-\frac34 \ln p} \approx 0.4704$ and $p \approx 0.7445$ is the unique positive root of $3p^6 + 2p - 2 = 0$. 
\end{prop}

We also obtained the following lower bound example for 5 colors through interactions with AI, but were unable to verify its optimality. At this time, we are not sure how to go about proving that this is optimal; solving the optimization program in \cref{conj:optimization} would be very hard, even with AI assistance. For this construction, we partition $V(G)$ into nine equal-size components that we identify with the $3 \times 3$ grid or affine plane $\mathbb{F}_3^2$. The five partitions $\mc{P}_1$, \ldots, $\mc{P}_5$ in the optimization problem are $\set{[9]}$ and the partitions into lines in four different directions: the ``row''  $\{\{1,2,3\}, \{4,5,6\}, \{7,8,9\}\}$, the ``column'' partition $\{\{1,4,7\}, \{2,5,8\}, \{3,6,9\}\}$, and the two ``diagonal'' partitions $\{\{1,5,9\}, \{2,6,7\}, \{3,4,8\}\}$ and $\{\{1,6,8\}, \{2,4,9\}, \{3,5,7\}\}$. Thus any pair of components lie in exactly one of the nontrivial color-specific partitions. Consequently, across components, an edge can only take one of two colors: the spanning color or the specific non-spanning color that has edges between those two parts, as was the case for four colors.

The best choice of parameters that AI could find was similar to the optimal four-color example. Within a component, an edge takes any of the four non-spanning colors with probability $\frac14$ each, then between any two components, an edge takes the corresponding non-spanning color with probability $p$ (and the spanning color with probability $1-p$). With this restriction, the lower bound appears to be optimized when $p \approx 0.8456$ is the unique positive root of $4p^{24} = 3(1-p)^2$, yielding $\beta_5 \leq 0.3861$.

We finish with the table of values of $\beta_{\ell}^{-1}$ such that	$R_{\had}(k;\ell)=(\beta_\ell^{-1}+o(1))k\sqrt{\ln k}$. The third column shows the lower bound obtained by considering a uniformly random coloring of edges into $\ell$ colors, which is optimal for two colors as shown in \cite{hadramtwo}.

\begin{table}[htbp]
	\centering
	\renewcommand{\arraystretch}{1.2}
	\begin{tabular}{c c c l}
		\hline
		$\ell$ & $\beta_\ell^{-1}$
		& \shortstack{uniform lower bound}
		& source \\  \hline
		$2$ & $1.2011\ldots$      & $1.2011\ldots$ & \cite{hadramtwo} \\
		$3$ & $1.5814\ldots$      & $1.5704\ldots$ & Theorem~\ref{t:main} \\
		$4$ & $2.1259\ldots$      & $1.8644\ldots$ & Proposition~\ref{conj:fourcol} \\
		$5$ & $\geq 2.5901\ldots$ & $2.1169\ldots$ & Example above  \\
		\hline
	\end{tabular}
\end{table}

\printbibliography

\newpage

\appendix
\counterwithin{thm}{section}
\renewcommand{\thethm}{\Alph{section}.\arabic{thm}}

\section{Proof of \cref{thomasonupper}}

Let $G$ be a graph on $n$ vertices with a partition $(X_1,\dots,X_k)$ of $V(G)$ into sets of size at least two, and let $\alpha_i$ and $p_{ij}$ be as in \cref{thomasonupper}. We write $q_{ij}=1-p_{ij}$ for the density of nonadjacent pairs, and for $0<c\leq 1$ we let
\[
 \Lambda_c(G)=-\sum_{(i,j)\in[k]^2}\alpha_i\alpha_j\log\max\{c,q_{ij}\}.
\]
The partition will always be clear from the context; for $Z\subseteq V(G)$, the graph $G-Z$ is equipped with the partition $(X_1\setminus Z,\dots,X_k\setminus Z)$. Since $1-p^*_{ij}=\max\{\eps,q_{ij}\}$, \cref{thomasonupper} asserts that
\[
 \had(G)\geq\frac{n}{\sqrt{\log n}}\left(\sqrt{\Lambda_\eps(G)}-\eps\right)
\]
whenever $n\geq n_0(\eps)$, $\kappa(G)\geq\eps n$, and $|X_i|\geq\eps n$ for every $i$.

A partition of $V(G)$ is \defn{complete} if $G$ has an edge between any two of its parts, and it is an \defn{equipartition} if its part sizes differ by at most one. The proof follows Thomason's complete-equipartition and dense-minor arguments \cite[Theorem 3.1, Corollary 3.1, and Theorem 4.1]{MR1814910}. After removing two small random sets of vertices, we find a complete equipartition of the remaining graph into nearly $\sqrt{\Lambda_\eps(G)}\,n/\sqrt{\log n}$ parts of size $O(\sqrt{\log n})$. We then use the removed vertices and the connectivity of $G$ to extend these parts to the branch sets of a complete minor.

\subsection{Complete equipartitions}

The next lemma is a multitype analogue of Thomason's complete-equipartition theorem: every part of a random partition receives a prescribed number $\ell_i$ of vertices from $X_i$. 

\begin{lem}\label[lem]{lem:multitype-finite}
Let $s\geq2$ and $\ell_1,\dots,\ell_k\geq1$ be integers, and let $\ell=\ell_1+\dots+\ell_k$. Let $G$ be a graph with a partition $(X_1,\dots,X_k)$ of $V(G)$ such that $|X_i|=s\ell_i$ for every $i$. Let $0<\tau<1$, $0<\varphi<1$, and $\omega\geq1$, and let
\begin{equation}\label{eq:P-def}
 P=(2\omega)^\ell\,\tau^{-\ell(k+\varphi\ell)}\prod_{(i,j)\in[k]^2}\min\{1,q_{ij}+\tau\}^{\ell_i\ell_j}.
\end{equation}
Then $G$ has a complete equipartition into at most $s$ and at least
\[
 s-\frac{4s}{\omega\varphi}-4s^2P
\]
parts.
\end{lem}

\begin{proof}
If $s-4s/(\omega\varphi)-4s^2P\leq1$, the partition with a single part satisfies the lemma. We may therefore assume otherwise.

We construct a random partition of $V(G)$ into $s$ parts of size $\ell$ and show that, in expectation, few pairs of parts have no edge between them. For a vertex $u$, let $Q(u)$ be the set of vertices $v\neq u$ not adjacent to $u$, and let
\[
 \rho(u)=\prod_{j=1}^k\min\left\{1,\frac{|Q(u)\cap X_j|}{|X_j|}+\tau\right\}^{\ell_j/\ell}.
\]
Then $\tau\leq\rho(u)\leq1$. Fix $i$ and $j$. Counting nonadjacent pairs gives $\sum_{u\in X_i}|Q(u)\cap X_j|=q_{ij}|X_i||X_j|$ if $i\neq j$, and $\sum_{u\in X_i}|Q(u)\cap X_i|=q_{ii}|X_i|(|X_i|-1)$. Thus the average of $|Q(u)\cap X_j|/|X_j|$ over $u\in X_i$ is at most $q_{ij}$. As the average of a minimum is at most the minimum of the averages, the average of $\min\{1,|Q(u)\cap X_j|/|X_j|+\tau\}$ over $u\in X_i$ is at most $\min\{1,q_{ij}+\tau\}$. The inequality of arithmetic and geometric means therefore gives
\begin{equation}\label{eq:rho-geomean}
 \prod_{u\in X_i}\rho(u)\leq\prod_{j=1}^k\min\{1,q_{ij}+\tau\}^{s\ell_i\ell_j/\ell}
 \qquad(1\leq i\leq k).
\end{equation}

Order each class $X_i$ so that $\rho$ is nondecreasing, and divide it into consecutive blocks $B_{i,1},\dots,B_{i,\ell_i}$ of $s$ vertices each. Choose, independently for every block, a uniformly random ordering of its vertices, and for $1\leq t\leq s$ let $W_t$ consist of the $t$th vertices of all blocks. Thus $(W_1,\dots,W_s)$ is a partition of $V(G)$ into $s$ sets of size $\ell$. For $W\subseteq V(G)$, let $N(W)$ be the set of vertices outside $W$ with no neighbor in $W$; thus $u\in N(W)$ if and only if $W\subseteq Q(u)$.

Fix a vertex $u$ and an index $t$. The vertex of $W_t$ in a block $B_{j,b}$ is uniformly distributed in that block, independently for different blocks. By the inequality of arithmetic and geometric means applied to the blocks of each class,
\[
 \Prob(W_t\subseteq Q(u))
 =\prod_{j=1}^k\prod_{b=1}^{\ell_j}\frac{|Q(u)\cap B_{j,b}|}{s}
 \leq\prod_{j=1}^k\left(\frac{|Q(u)\cap X_j|}{s\ell_j}\right)^{\ell_j}
 \leq\rho(u)^\ell.
\]
Let $\rho_{i,b}=\max_{u\in B_{i,b}}\rho(u)$. Summing over $u\in B_{i,b}$ gives
\begin{equation}\label{eq:block-expectation}
 \mathbb E\,|B_{i,b}\cap N(W_t)|\leq s\rho_{i,b}^{\ell}.
\end{equation}

We say that $W_t$ \defn{rejects} a block $B_{i,b}$ with $b<\ell_i$ if $|B_{i,b}\cap N(W_t)|>\omega s\rho_{i,b}^{\ell}$, and that $W_t$ is \defn{acceptable} if it rejects at most $\varphi\ell$ blocks. By \eqref{eq:block-expectation} and Markov's inequality, each block is rejected with probability at most $1/\omega$. There are fewer than $\ell$ blocks $B_{i,b}$ with $b<\ell_i$, so the expected number of rejected blocks is less than $\ell/\omega$, and Markov's inequality gives
\begin{equation}\label{eq:unacceptable-prob}
 \Prob(W_t\text{ is not acceptable})\leq\frac1{\omega\varphi}.
\end{equation}

Fix an acceptable value $W$ of $W_t$, and for every $i$ let $M_i$ be the set of indices $b<\ell_i$ such that $W$ does not reject $B_{i,b}$. Since $W$ rejects at most $\varphi\ell$ of the $\ell-k$ blocks $B_{i,b}$ with $b<\ell_i$,
\begin{equation}\label{eq:M-lower}
 \sum_{i=1}^k\left(\ell_i-|M_i|\right)\leq k+\varphi\ell.
\end{equation}
For $b<\ell_i$, every vertex of $B_{i,b+1}$ has $\rho$-value at least $\rho_{i,b}$ by the choice of the ordering of $X_i$, and so $\rho_{i,b}^s\leq\prod_{u\in B_{i,b+1}}\rho(u)$. Multiplying over $b\in M_i$ and using $\rho\geq\tau$ on the remaining $\ell_i-|M_i|$ blocks of $X_i$, we obtain
\[
 \Bigl(\prod_{b\in M_i}\rho_{i,b}\Bigr)^s\leq\tau^{-s(\ell_i-|M_i|)}\prod_{u\in X_i}\rho(u).
\]
We take $s$th roots, apply \eqref{eq:rho-geomean}, multiply over $i$, and raise to the power $\ell$. By \eqref{eq:M-lower} and $\tau<1$, this gives
\begin{equation}\label{eq:block-product}
 \prod_{i=1}^k\prod_{b\in M_i}\rho_{i,b}^{\ell}
 \leq\tau^{-\ell(k+\varphi\ell)}\prod_{(i,j)\in[k]^2}\min\{1,q_{ij}+\tau\}^{\ell_i\ell_j}.
\end{equation}

Let $t'\neq t$. Conditional on $W_t=W$, the vertex of $W_{t'}$ in each block is uniformly distributed among the $s-1$ vertices of the block outside $W$, independently for different blocks. Since $W_{t'}$ is disjoint from $W$, there is no edge between $W_{t'}$ and $W$ if and only if $W_{t'}\subseteq N(W)$. As $N(W)\cap W=\varnothing$, considering only the blocks $B_{i,b}$ with $b\in M_i$ gives
\[
 \Prob\bigl(W_{t'}\subseteq N(W)\bigm|W_t=W\bigr)
 \leq\prod_{i=1}^k\prod_{b\in M_i}\frac{|B_{i,b}\cap N(W)|}{s-1}
 \leq\prod_{i=1}^k\prod_{b\in M_i}\frac{\omega s}{s-1}\rho_{i,b}^{\ell}
 \leq P,
\]
where the last inequality follows from $s/(s-1)\leq2$, $\sum_i|M_i|\leq\ell$, \eqref{eq:block-product}, and \eqref{eq:P-def}.

Consequently, for $t\neq t'$, the probability that $W_t$ is acceptable and there is no edge between $W_t$ and $W_{t'}$ is at most $P$. By \eqref{eq:unacceptable-prob}, the expected number of parts that are not acceptable is at most $s/(\omega\varphi)$, and the expected number of pairs of acceptable parts with no edge between them is at most $\binom s2P<s^2P$. By Markov's inequality, there is an outcome with at most $4s/(\omega\varphi)$ parts that are not acceptable and at most $4s^2P$ pairs of acceptable parts with no edge between them. We discard all parts that are not acceptable and one part from each such pair. The at least $s-4s/(\omega\varphi)-4s^2P$ remaining parts are pairwise joined by edges. Distributing the vertices of the discarded parts among them, one at a time and always to a smallest part, gives a complete equipartition with the required number of parts.
\end{proof}

The next lemma shows that deleting few vertices changes $\Lambda_a$ only slightly. We apply it to pass to class sizes divisible by prescribed integers, and to remove the random sets of vertices in the proof of \cref{thomasonupper}.

\begin{lem}\label[lem]{lem:stability}
Let $0<a\leq1$ and $0\leq\delta\leq a/4$. Let $G$ be a graph on $n\geq4/a$ vertices with a partition $(X_1,\dots,X_k)$ of $V(G)$ such that $|X_i|\geq an$ for every $i$, and let $Z\subseteq V(G)$ with $|Z|\leq\delta n$. Then
\[
 \left|\Lambda_a(G-Z)-\Lambda_a(G)\right|\leq\frac{12\delta}{a^3}.
\]
\end{lem}

\begin{proof}
Primes denote the quantities defined for $G-Z$. Each class of $G-Z$ has at least $(a-\delta)n\geq3an/4$ vertices.

We first compare the densities. Fix $i,j\in[k]$. Let $N$ and $N'$ be the numbers of pairs of distinct vertices with one vertex in $X_i$ and the other in $X_j$, in $G$ and in $G-Z$, respectively, and let $e$ and $e'$ be the numbers of such pairs that are edges. Every pair counted by $N$ but not by $N'$ contains a vertex of $Z$, so $0\leq e-e'\leq N-N'\leq|Z|n\leq\delta n^2$. Since $e\leq N$, both $(e-e')N$ and $e(N-N')$ lie in $[0,N(N-N')]$, and so
\[
 \left|\frac{e'}{N'}-\frac eN\right|=\frac{\left|(e-e')N-e(N-N')\right|}{NN'}\leq\frac{N-N'}{N'}.
\]
In both cases $i\neq j$ and $i=j$, as $an\geq4$, we have $N'\geq\frac{3an}{4}\left(\frac{3an}{4}-1\right)/2\geq\frac{3a^2n^2}{16}$. Therefore
\[
 |q'_{ij}-q_{ij}|=|p'_{ij}-p_{ij}|\leq\frac{16\delta}{3a^2}\leq\frac{6\delta}{a^2}.
\]

Next, $\alpha'_i-\alpha_i=\bigl(|X_i||Z|-n|Z\cap X_i|\bigr)/\bigl(n(n-|Z|)\bigr)$, and so
\[
 \sum_{i=1}^k|\alpha'_i-\alpha_i|\leq\frac{2n|Z|}{n(n-|Z|)}\leq\frac{2\delta}{1-\delta}\leq3\delta.
\]
Consequently,
\[
 \sum_{(i,j)\in[k]^2}|\alpha'_i\alpha'_j-\alpha_i\alpha_j|
 \leq\sum_{(i,j)\in[k]^2}\left(|\alpha'_i-\alpha_i|\,\alpha'_j+\alpha_i\,|\alpha'_j-\alpha_j|\right)
 =2\sum_{i=1}^k|\alpha'_i-\alpha_i|\leq6\delta.
\]

Finally, the function $f(x)=-\log\max\{a,x\}$ is $(1/a)$-Lipschitz on $[0,1]$ and satisfies $0\leq f\leq\log(1/a)\leq1/a$ there. Hence
\[
 \left|\Lambda_a(G-Z)-\Lambda_a(G)\right|
 \leq\sum_{(i,j)\in[k]^2}\alpha'_i\alpha'_j\,|f(q'_{ij})-f(q_{ij})|
 +\sum_{(i,j)\in[k]^2}|\alpha'_i\alpha'_j-\alpha_i\alpha_j|\,f(q_{ij})
 \leq\frac{6\delta}{a^3}+\frac{6\delta}{a}\leq\frac{12\delta}{a^3}.\qedhere
\]
\end{proof}

Combining the two lemmas, we obtain the complete equipartition used in the proof of \cref{thomasonupper}.

\begin{lem}\label[lem]{lem:multitype-asymptotic}
For all $0<a<1$, $\lambda>0$, and $0<\eta<1$ there exists $N$ such that the following holds. Let $G$ be a graph on $n\geq N$ vertices with a partition $(X_1,\dots,X_k)$ of $V(G)$ such that $|X_i|\geq an$ for every $i$, and suppose that $\Lambda=\Lambda_a(G)\geq\lambda$. Then $G$ has a complete equipartition into $m$ parts, where
\[
 (1-\eta)\,n\sqrt{\frac{\Lambda}{\log n}}\leq m\leq n\sqrt{\frac{\Lambda}{\log n}}.
\]
\end{lem}

\begin{proof}
We have $k\leq1/a$ and $\Lambda\leq\log(1/a)$. Let $\xi=\min\{1/20,\eta/7\}$ and $\tau=\min\{1/2,a\xi\lambda\}$, choose $0<\varphi<1$ with $\varphi\log(1/\tau)\leq\xi\lambda$, and choose $\omega\geq1$ with $4/(\omega\varphi)\leq\xi$. These choices depend only on $a$, $\lambda$, and $\eta$, and we assume throughout that $n$ is sufficiently large in terms of $a$, $\lambda$, and $\eta$.

Let
\[
 \ell=\left\lceil(1+3\xi)\sqrt{\frac{\log n}{\Lambda}}\,\right\rceil.
\]
Since $\lambda\leq\Lambda\leq\log(1/a)$, we have $\ell\geq\sqrt{\log n/\log(1/a)}$ and $\ell\leq2\sqrt{\log n/\lambda}$. Choose integers $\ell_1,\dots,\ell_k$ with $\sum_i\ell_i=\ell$ and $|\ell_i-\alpha_i\ell|\leq1$ for every $i$, by rounding every $\alpha_i\ell$ down and then increasing some of the resulting integers by one. Then $\ell_i\geq\lfloor a\ell\rfloor\geq1$. Let $s=\min_i\lfloor|X_i|/\ell_i\rfloor$. Since $\ell_i\leq\alpha_i\ell+1\leq\alpha_i\ell(1+1/(a\ell))$, we have $|X_i|/\ell_i\geq(n/\ell)(1-1/(a\ell))$ for every $i$. As $n/(a\ell^2)\geq1$, it follows that
\begin{equation}\label{eq:s-bounds}
 \frac n\ell\left(1-\frac2{a\ell}\right)\leq s\leq\frac n\ell,
\end{equation}
where the upper bound holds because $s\ell=\sum_is\ell_i\leq n$.

Choose $X_i^0\subseteq X_i$ with $|X_i^0|=s\ell_i$, and let $G_0=G[X_1^0\cup\dots\cup X_k^0]$. Then $G_0=G-Z$ with $|Z|=n-s\ell\leq2n/(a\ell)$ by \eqref{eq:s-bounds}, and the class $X_i^0$ of $G_0$ has proportion $\ell_i/\ell$. By \cref{lem:stability} with $\delta=2/(a\ell)$,
\[
 \Lambda_a(G_0)\geq\Lambda-\frac{24}{a^4\ell}\geq\Lambda-\xi\lambda\geq(1-\xi)\Lambda.
\]

We apply \cref{lem:multitype-finite} to $G_0$ with the integers $s,\ell_1,\dots,\ell_k$ and the parameters $\tau,\varphi,\omega$. Let $q^0_{ij}$ be the densities of nonadjacent pairs in $G_0$. For $0\leq x\leq1$ we have $\min\{1,x+\tau\}\leq\max\{a,x\}+\tau\leq\max\{a,x\}(1+\tau/a)$, and so $\log\min\{1,x+\tau\}\leq\log\max\{a,x\}+\tau/a$. As the classes of $G_0$ have proportions $\ell_i/\ell$ and $\tau/a\leq\xi\lambda\leq\xi\Lambda$,
\[
 \prod_{(i,j)\in[k]^2}\min\{1,q^0_{ij}+\tau\}^{\ell_i\ell_j}
 \leq\exp\left(-\ell^2\left(\Lambda_a(G_0)-\frac\tau a\right)\right)
 \leq\exp\left(-(1-2\xi)\Lambda\ell^2\right).
\]
Using $k\leq1/a$ and $\varphi\log(1/\tau)\leq\xi\Lambda$, we conclude that the quantity $P$ in \eqref{eq:P-def} satisfies
\[
 \log P\leq\ell\log(2\omega)+\ell(k+\varphi\ell)\log\frac1\tau-(1-2\xi)\Lambda\ell^2
 \leq\ell\left(\log(2\omega)+\frac1a\log\frac1\tau\right)-(1-3\xi)\Lambda\ell^2.
\]
By the choice of $\ell$, we have $\Lambda\ell^2\geq(1+3\xi)^2\log n$, and $(1-3\xi)(1+3\xi)^2\geq1+2\xi$ as $\xi\leq1/20$. Since $\ell\leq2\sqrt{\log n/\lambda}$, it follows that $\log P\leq-(1+\xi)\log n$, and so $4s^2P\leq4sn^{-\xi}\leq\xi s$. Together with $4s/(\omega\varphi)\leq\xi s$, \cref{lem:multitype-finite} gives a complete equipartition of $G_0$ into $m$ parts with $(1-2\xi)s\leq m\leq s$. Adding the vertices of $Z$ one at a time, each to a smallest part, gives a complete equipartition of $G$ into $m$ parts.

It remains to bound $m$. On the one hand, $m\leq s\leq n/\ell\leq n\sqrt{\Lambda/\log n}$. On the other hand, $\sqrt{\log n/\Lambda}\geq\sqrt{\log n/\log(1/a)}$ tends to infinity, so $\ell\leq(1+4\xi)\sqrt{\log n/\Lambda}$ and $2/(a\ell)\leq\xi$ for large $n$. By \eqref{eq:s-bounds},
\[
 m\geq(1-2\xi)(1-\xi)\frac n\ell\geq\frac{1-3\xi}{1+4\xi}\,n\sqrt{\frac\Lambda{\log n}}\geq(1-7\xi)\,n\sqrt{\frac\Lambda{\log n}}\geq(1-\eta)\,n\sqrt{\frac\Lambda{\log n}}.\qedhere
\]
\end{proof}

\subsection{Connecting the parts}

We use the following Chernoff bounds \cite[Theorems 4.4 and 4.5]{Mitzenmacher_Upfal_2005}. If $Y$ is a sum of independent Bernoulli random variables, then
\begin{equation}\label{eq:Chernoff}
 \Prob(Y\geq2\mu)\leq e^{-\mu/3}\quad\text{for }\mu=\mathbb EY,
 \qquad
 \Prob(Y\leq\mu/2)\leq e^{-\mu/8}\quad\text{for }0<\mu\leq\mathbb EY.
\end{equation}
For $\mu<\mathbb EY$, the second bound follows from the case $\mu=\mathbb EY$.

The two random sets of vertices are provided by the following lemma. For each part, $D$ supplies a small set containing a neighbor of every vertex of the part, and $C$ supplies short paths joining the vertices of this set.

\begin{lem}\label[lem]{lem:reservoirs}
Let $0<c<1$ and $h=\lceil2/c\rceil+1$. For all sufficiently large $n$, with $r=1/\log\log n$, every $n$-vertex graph $G$ with $\kappa(G)\geq cn$ contains disjoint sets $C,D\subseteq V(G)$ such that
\begin{enumerate}[(i)]
 \item $|C|,|D|\leq2rn$;
 \item every vertex of $G$ has at least $crn/2$ neighbors in $D$;
 \item any two distinct vertices of $G$ are joined by at least $nr^h$ internally disjoint paths of length at most $h$ whose internal vertices all lie in $C$.
\end{enumerate}
\end{lem}

\begin{proof}
Independently for every vertex, place it in $C$ with probability $r$, in $D$ with probability $r$, and in neither set otherwise. By \eqref{eq:Chernoff}, (i) fails with probability at most $2e^{-rn/3}$.

Every vertex $v$ has degree at least $\kappa(G)\geq cn$, so the expected number of neighbors of $v$ in $D$ is at least $crn$. By \eqref{eq:Chernoff}, $v$ has fewer than $crn/2$ neighbors in $D$ with probability at most $e^{-crn/8}$. Thus (ii) fails with probability at most $ne^{-crn/8}$.

Fix distinct vertices $u$ and $v$. By Menger's theorem (see, for example, \cite[Section 3.3]{MR3644391}), $G$ contains $\kappa(G)$ internally disjoint $u$--$v$ paths. Their sets of internal vertices are disjoint and contain fewer than $n$ vertices in total, so fewer than $\kappa(G)/2$ of these paths have more than $2n/\kappa(G)$ internal vertices. Hence more than $\kappa(G)/2\geq cn/2$ of them have at most $2n/\kappa(G)\leq 2/c$, and hence at most $h-1$, internal vertices. Each of these paths has all its internal vertices in $C$ with probability at least $r^{h-1}$, independently of the others. The expected number of such paths is at least $cnr^{h-1}/2\geq2nr^h$, as $4r\leq c$ for large $n$. By \eqref{eq:Chernoff}, fewer than $nr^h$ of these paths have all their internal vertices in $C$ with probability at most $\exp(-cnr^{h-1}/16)$. Thus (iii) fails with probability at most $n^2\exp(-cnr^{h-1}/16)$.

Since $rn$ and $nr^{h-1}=n/(\log\log n)^{h-1}$ are both much larger than $\log n$, the three failure probabilities sum to less than one for large $n$.
\end{proof}

\begin{proof}[Proof of \cref{thomasonupper}]
If $\eps\geq1$, no graph on $n\geq2$ vertices satisfies $\kappa(G)\geq\eps n$, so we assume $0<\eps<1$. Let $G$ and $(X_1,\dots,X_k)$ satisfy the hypotheses, and let $\Lambda^*=\Lambda_\eps(G)$. We must show that
\begin{equation}\label{eq:upper-goal}
 \had(G)\geq\frac n{\sqrt{\log n}}\left(\sqrt{\Lambda^*}-\eps\right).
\end{equation}
Since $\max\{\eps,q_{ij}\}\geq\eps$, we have $\Lambda^*\leq\log(1/\eps)$. If $\Lambda^*\leq\eps^2$, the right-hand side of \eqref{eq:upper-goal} is not positive, so we assume $\Lambda^*>\eps^2$. Let
\[
 \eta=\min\left\{\frac12,\frac{\eps}{2\sqrt{\log(1/\eps)}}\right\},
 \qquad
 h=\left\lceil\frac2\eps\right\rceil+1,
 \qquad
 r=\frac1{\log\log n},
\]
and assume that $n$ is sufficiently large in terms of $\eps$.

\emph{The complete equipartition.} Let $C$ and $D$ be the sets given by \cref{lem:reservoirs} with $c=\eps$, let $G'=G-(C\cup D)$, and let $n'=|V(G')|$, so that $n'\geq(1-4r)n$. Every class of $G'$ has at least $\eps n-4rn\geq\eps n'/2$ vertices. By \cref{lem:stability} with $a=\eps$ and $\delta=4r$,
\[
 \Lambda_\eps(G')\geq\Lambda^*-\frac{48r}{\eps^3};
\]
in particular, $\Lambda_\eps(G')\geq\eps^2/2$. As $\max\{\eps/2,x\}\leq\max\{\eps,x\}$, we have $\eps^2/2\leq\Lambda_\eps(G')\leq\Lambda_{\eps/2}(G')\leq\log(2/\eps)$. By \cref{lem:multitype-asymptotic} with $a=\eps/2$, $\lambda=\eps^2/2$, and $\eta$ as above, $V(G')$ has a complete equipartition $(W_1,\dots,W_m)$ with
\begin{equation}\label{eq:m-bounds}
 (1-\eta)\,n'\sqrt{\frac{\Lambda_\eps(G')}{\log n'}}\leq m\leq n'\sqrt{\frac{\log(2/\eps)}{\log n'}}.
\end{equation}
Since $\Lambda^*>\eps^2$, the lower bounds on $n'$ and $\Lambda_\eps(G')$ give $n'\sqrt{\Lambda_\eps(G')}\geq(1-4r)\sqrt{1-48r/\eps^5}\cdot n\sqrt{\Lambda^*}$. Thus, for large $n$, using $\log n'\leq\log n$, $\sqrt{\Lambda^*}\leq\sqrt{\log(1/\eps)}$, and $2\eta\sqrt{\log(1/\eps)}\leq\eps$,
\[
 m\geq(1-\eta)^2\,n\sqrt{\frac{\Lambda^*}{\log n}}\geq(1-2\eta)\frac{n\sqrt{\Lambda^*}}{\sqrt{\log n}}\geq\frac n{\sqrt{\log n}}\left(\sqrt{\Lambda^*}-\eps\right).
\]
It therefore suffices to show that $G$ has a $K_m$ minor.

We also record two upper bounds. Since $n'\geq n/2$, we have $\log n'\geq(\log n)/2$, and the upper bound in \eqref{eq:m-bounds} gives
\begin{equation}\label{eq:m-upper}
 m\leq2\sqrt{\log(2/\eps)}\,\frac n{\sqrt{\log n}}.
\end{equation}
The lower bound in \eqref{eq:m-bounds}, together with $\eta\leq1/2$ and $\Lambda_\eps(G')\geq\eps^2/2$, gives $n'/m\leq2\sqrt2\,\eps^{-1}\sqrt{\log n'}$. Every part of an equipartition has at most $n'/m+1$ vertices, and so
\begin{equation}\label{eq:part-size}
 |W_b|\leq\frac4\eps\sqrt{\log n}\qquad(1\leq b\leq m).
\end{equation}

\emph{Dominating sets.} We now choose pairwise disjoint sets $M_1,\dots,M_m\subseteq D$ such that every vertex of $W_b$ has a neighbor in $M_b$ and
\begin{equation}\label{eq:M-size}
 |M_b|\leq\frac8\eps\log\log n\qquad(1\leq b\leq m).
\end{equation}
By \eqref{eq:m-upper}, sets satisfying \eqref{eq:M-size} have total size at most
\begin{equation}\label{eq:M-total}
 \frac{8m}{\eps}\log\log n\leq\frac{16\sqrt{\log(2/\eps)}}{\eps}\cdot\frac{n\log\log n}{\sqrt{\log n}}.
\end{equation}
For large $n$, the right-hand side of \eqref{eq:M-total} is at most $\eps rn/4=\eps n/(4\log\log n)$, and $h-1$ times it is less than $nr^h=n/(\log\log n)^h$.

Suppose that $M_1,\dots,M_{b-1}$ have been chosen, and let $U$ be their union. Then $|U|\leq\eps rn/4$, so by \cref{lem:reservoirs}(ii) and (i), every vertex of $W_b$ has at least $\eps rn/4\geq(\eps/8)|D\setminus U|$ neighbors in $D\setminus U$. We choose vertices of $D\setminus U$ one at a time, each time choosing a vertex adjacent to as many vertices of $W_b$ without a neighbor among the previously chosen vertices as possible, and we stop when no such vertex of $W_b$ remains. If $A\subseteq W_b$ is the set of such vertices before a step, then there are at least $(\eps/8)|A||D\setminus U|$ edges between $A$ and $D\setminus U$, and so the chosen vertex has at least $(\eps/8)|A|$ neighbors in $A$. Hence, for $j\geq1$, after $j$ steps fewer than $e^{-\eps j/8}|W_b|$ vertices of $W_b$ have no neighbor among the chosen vertices, and so the procedure stops after at most $\lfloor(8/\eps)\log|W_b|\rfloor+1$ steps. Let $M_b$ be the set of chosen vertices. By \eqref{eq:part-size},
\[
 |M_b|\leq1+\frac8\eps\log\left(\frac4\eps\sqrt{\log n}\right)\leq\frac8\eps\log\log n
\]
for large $n$. This completes the construction.

\emph{Connecting paths.} For every $b$, enumerate $M_b=\{x_1,\dots,x_{|M_b|}\}$. For all $b$ and $1\leq j<|M_b|$, we choose an $x_j$--$x_{j+1}$ path of length at most $h$ whose internal vertices lie in $C$, so that the sets of internal vertices of all chosen paths are pairwise disjoint. We choose the paths one at a time. Fewer than $\sum_b|M_b|$ paths are needed, each with at most $h-1$ internal vertices, so by \eqref{eq:M-total} fewer than $nr^h$ vertices of $C$ are used in total. When a path is to be chosen, \cref{lem:reservoirs}(iii) provides at least $nr^h$ internally disjoint candidates. Every previously used vertex is an internal vertex of at most one candidate, so some candidate avoids all previously used vertices; we choose it. Let $N_b$ be the set of internal vertices of the paths chosen for $M_b$.

\emph{The branch sets.} The sets $W_b\cup M_b\cup N_b$, for $1\leq b\leq m$, are pairwise disjoint: the sets $W_b$ partition $V(G)\setminus(C\cup D)$, the sets $M_b$ are disjoint subsets of $D$, and the sets $N_b$ are disjoint subsets of $C$. Each graph $G[M_b\cup N_b]$ is connected, and every vertex of $W_b$ has a neighbor in $M_b$, so each graph $G[W_b\cup M_b\cup N_b]$ is connected. Since $(W_1,\dots,W_m)$ is a complete partition of $V(G')$ and $G'\subseteq G$, there is an edge of $G$ between any two of these sets. Hence they are the branch sets of a $K_m$ minor of $G$, and $\had(G)\geq m$.
\end{proof}

\section{Proof of \cref{thomasonlower}}

Throughout this appendix, $G$ is the random graph of \cref{thomasonlower}, with vertex set $V=[n]$. For distinct $u,v\in V$, let $p_{uv}$ be the probability that $uv\in E(G)$; thus $p_{uv}=p_{ij}$ if $u\in X_i$ and $v\in X_j$. As in the proof of \cref{thomasonupper}, a partition of $V$ is \defn{complete} if $G$ has an edge between any two of its parts.

We adapt the argument of Bollob\'as, Catlin and Erd\H{o}s~\cite{MR593989} for $G(n,p)$ to unequal edge probabilities. Their calculation treats $p=\tfrac12$, and the paper also records the extension to every fixed $0<p<1$ \cite[p.~198]{MR593989}. The following estimate involves no asymptotic assumptions.

\begin{lem}\label[lem]{lem:finite-estimate}
For every integer $s\geq2$,
\[
 \Prob\bigl(\had(G)\geq s\bigr)\leq s^n\exp\left(-\binom s2\exp\left(-\frac{n^2L^2}{s(s-1)}\right)\right).
\]
\end{lem}

\begin{proof}
If $s>n$, the left-hand side is zero, so we assume $s\leq n$. If $G$ has a $K_s$ minor, then adding all vertices outside its branch sets to one of the branch sets gives a complete partition of $V$ into $s$ nonempty parts. There are at most $s^n$ ordered partitions of $V$ into $s$ parts. By the union bound, it therefore suffices to show that a fixed ordered partition $(A_1,\dots,A_s)$ of $V$ into nonempty sets is complete with probability at most $\exp\bigl(-\binom s2\exp(-n^2L^2/(s(s-1)))\bigr)$.

For $1\leq a<b\leq s$, let
\[
 W_{ab}=\sum_{u\in A_a}\sum_{v\in A_b}-\log(1-p_{uv}).
\]
The probability that $G$ has no edge between $A_a$ and $A_b$ is $\prod_{u\in A_a}\prod_{v\in A_b}(1-p_{uv})=e^{-W_{ab}}$. These events are independent for distinct pairs $\{a,b\}$, as they involve disjoint sets of pairs of vertices. Hence, by the inequality $1-x\leq e^{-x}$,
\begin{equation}\label{eq:complete-prob}
 \Prob\bigl((A_1,\dots,A_s)\text{ is complete}\bigr)=\prod_{1\leq a<b\leq s}\left(1-e^{-W_{ab}}\right)\leq\exp\left(-\sum_{1\leq a<b\leq s}e^{-W_{ab}}\right).
\end{equation}

Every pair of distinct vertices contributes to at most one of the sums $W_{ab}$, and so
\[
 \sum_{1\leq a<b\leq s}W_{ab}\leq\frac12\sum_{(i,j)\in[k]^2}|X_i||X_j|\bigl(-\log(1-p_{ij})\bigr)=\frac{n^2L^2}{2};
\]
indeed, the middle sum contains the term $-\log(1-p_{uv})$ for every ordered pair $(u,v)$ of distinct vertices, together with nonnegative terms for the pairs $(u,u)$.
As $x\mapsto e^{-x}$ is convex, Jensen's inequality gives
\[
 \sum_{1\leq a<b\leq s}e^{-W_{ab}}\geq\binom s2\exp\left(-\frac{1}{\binom s2}\sum_{1\leq a<b\leq s}W_{ab}\right)\geq\binom s2\exp\left(-\frac{n^2L^2}{s(s-1)}\right).
\]
Substituting this bound into \eqref{eq:complete-prob} completes the proof.
\end{proof}

\begin{proof}[Proof of \cref{thomasonlower}]
If $\eps>1$, no $p_{ij}$ satisfies $0\leq p_{ij}\leq1-\eps$, and if $\eps=1$, then all $p_{ij}$ are zero and $L=0$. In both cases the assertion is vacuous, so we assume $0<\eps<1$. Let
\[
 s=\left\lceil(1+\eps)\frac{nL}{\sqrt{\log n}}\right\rceil.
\]
If $\had(G)<s$, then $\had(G)\leq(1+\eps)nL/\sqrt{\log n}$. It therefore suffices to show that $\Prob(\had(G)\geq s)\leq\eps$ for all sufficiently large $n$, where the bound on $n$ depends only on $\eps$.

Since $0\leq p_{ij}\leq1-\eps$, every term $-\log(1-p_{ij})$ lies in $[0,\log(1/\eps)]$. As $\sum_{(i,j)\in[k]^2}\alpha_i\alpha_j=1$ and $L>\eps$, we obtain $\eps^2<L^2\leq\log(1/\eps)$. In particular, $s\geq(1+\eps)\eps n/\sqrt{\log n}\geq1+1/\eps$ for large $n$. We may also assume $s\leq n$, as otherwise $\Prob(\had(G)\geq s)=0$. As $s\geq1+1/\eps$, we have $s-1\geq s/(1+\eps)$, and so
\[
 s(s-1)\geq\frac{s^2}{1+\eps}\geq(1+\eps)\frac{n^2L^2}{\log n}.
\]
Thus $n^2L^2/(s(s-1))\leq\log n/(1+\eps)$ and $\binom s2\geq\eps^2n^2/(2\log n)$. By \cref{lem:finite-estimate} and $s\leq n$,
\[
 \Prob\bigl(\had(G)\geq s\bigr)\leq\exp\left(n\log n-\frac{\eps^2n^2}{2\log n}\,n^{-1/(1+\eps)}\right)=\exp\left(n\log n-\frac{\eps^2n^{1+\eps/(1+\eps)}}{2\log n}\right).
\]
The right-hand side depends only on $n$ and $\eps$, and it tends to zero as $n\to\infty$. Hence it is at most $\eps$ for all sufficiently large $n$.
\end{proof}

\end{document}

%% file: three-color.tikz
\begingroup

\definecolor{Gone}{RGB}{206,44,58}     
\definecolor{Gtwo}{RGB}{38,102,192}    
\definecolor{Gthree}{RGB}{24,150,96}   

\begin{tikzpicture}[
    x=1.6cm, y=1.6cm,
    line cap=round,
    eone/.style   = {draw=Gone,   line width=0.65pt, draw opacity=0.9},
    etwo/.style   = {draw=Gtwo,   line width=0.65pt, draw opacity=0.95},
    ethree/.style = {draw=Gthree, line width=0.65pt, draw opacity=0.95},
    cross/.style  = {line width=0.42pt, draw opacity=0.75},
    vertex/.style = {circle, inner sep=0pt, minimum size=4pt,
                     fill=black!80, draw=white, line width=0.55pt},
    part/.style   = {draw=black!40, line width=0.6pt},
  ]
  \def\n{8}               
  \def\hgt{1.45}          
  \def\rx{1.3} \def\ry{1.9}  
  \def\dist{5.0}          
  \def\loose{2.1}         
  \foreach \i in {1,...,\n} {
    \pgfmathsetmacro\vy{\hgt*(2*(\i-1)/(\n-1)-1)}
    \pgfmathsetmacro\vx{\rx*sqrt(1-(\vy/\ry)^2)}
    \coordinate (x\i) at (\vx,\vy);
    \coordinate (y\i) at (\dist-\vx,\vy);
  }
  \foreach \cx in {0,\dist} {
    \shade[inner color=Gone!14, outer color=Gone!4]
      (\cx,0) ellipse [x radius=\rx, y radius=\ry];
    \draw[part] (\cx,0) ellipse [x radius=\rx, y radius=\ry];
  }
  \foreach \i/\j/\col in {
    1/3/etwo, 3/7/etwo, 4/2/ethree, 1/7/ethree, 6/3/etwo, 2/8/etwo, 3/3/etwo, 7/3/ethree,
    4/6/etwo, 6/4/etwo, 1/6/etwo, 5/5/etwo, 1/1/ethree, 4/1/etwo, 2/3/ethree, 6/1/ethree,
    5/3/etwo, 3/8/etwo, 8/8/ethree, 3/5/ethree, 7/1/ethree, 2/7/ethree, 3/4/ethree, 8/2/etwo,
    4/8/ethree, 7/2/etwo, 4/3/ethree, 2/2/etwo, 8/3/ethree, 5/2/ethree, 8/4/etwo, 1/4/ethree,
    2/4/etwo, 5/1/etwo, 6/6/etwo, 6/8/ethree, 6/2/etwo, 4/4/ethree, 7/6/ethree, 7/5/ethree,
    5/7/etwo, 6/5/ethree, 4/5/etwo, 3/6/ethree, 5/4/ethree, 2/5/etwo, 4/7/etwo, 7/7/etwo,
    6/7/ethree, 8/5/ethree, 8/7/ethree, 2/1/ethree, 3/1/etwo, 7/8/etwo, 3/2/ethree, 2/6/ethree,
    5/8/ethree, 1/2/ethree, 8/6/etwo, 5/6/ethree, 1/8/etwo, 7/4/etwo, 8/1/etwo, 1/5/etwo}
    \draw[\col, cross] (x\i) -- (y\j);
  \foreach \i/\j/\col in {
    1/8/ethree, 1/7/eone, 2/8/eone, 2/7/eone, 3/8/eone, 1/6/etwo, 2/6/eone, 4/8/eone,
    1/5/eone, 3/7/etwo, 3/6/eone, 1/4/eone, 4/7/eone, 5/8/eone, 2/5/ethree, 1/3/eone,
    5/7/eone, 4/6/eone, 2/4/eone, 6/8/eone, 3/5/eone, 1/2/eone, 5/6/eone, 4/5/eone,
    2/3/eone, 6/7/eone, 3/4/etwo, 7/8/ethree}
    \draw[\col] (x\i) to[out=180, in=180, looseness=\loose] (x\j);
  \foreach \i/\j/\col in {
    1/8/eone, 2/8/etwo, 1/7/ethree, 3/8/eone, 2/7/eone, 1/6/eone, 1/5/eone, 4/8/eone,
    2/6/etwo, 3/7/ethree, 4/7/eone, 2/5/eone, 5/8/eone, 3/6/eone, 1/4/eone, 4/6/eone,
    5/7/eone, 1/3/eone, 3/5/eone, 2/4/eone, 6/8/eone, 7/8/eone, 1/2/eone, 2/3/eone,
    3/4/eone, 6/7/eone, 5/6/etwo, 4/5/ethree}
    \draw[\col] (y\i) to[out=0, in=0, looseness=\loose] (y\j);
  \foreach \i in {1,...,\n} { \node[vertex] at (x\i) {}; \node[vertex] at (y\i) {}; }
  \node[below=3pt] at (0,-\ry)     {$X$};
  \node[below=3pt] at (\dist,-\ry) {$Y$};
\end{tikzpicture}
\endgroup

%% file: four-color.tikz
\begingroup

\definecolor{Gone}{RGB}{206,44,58}     
\definecolor{Gtwo}{RGB}{38,102,192}    
\definecolor{Gthree}{RGB}{24,150,96}   
\definecolor{Gfour}{RGB}{226,150,20}   

\begin{tikzpicture}[
    x=1.12cm, y=1.12cm,
    line cap=round,
    eone/.style   = {draw=Gone,   line width=0.65pt, draw opacity=0.9},
    etwo/.style   = {draw=Gtwo,   line width=0.65pt, draw opacity=0.95},
    ethree/.style = {draw=Gthree, line width=0.65pt, draw opacity=0.95},
    efour/.style  = {draw=Gfour,  line width=0.65pt, draw opacity=0.95},
    cross/.style  = {line width=0.36pt, draw opacity=0.72},
    vertex/.style = {circle, inner sep=0pt, minimum size=4pt,
                     fill=black!80, draw=white, line width=0.55pt},
    part/.style   = {draw=black!40, line width=0.6pt},
  ]
  \def\n{6}               
  \def\hgt{1.1}          
  \def\rx{1.05} \def\ry{1.55}  
  \def\off{2.8}           
  \def\loose{2.25}         
  \def\rotone{-45} \def\rottwo{-135} \def\rotthree{45} \def\rotfour{135}
  \foreach \q/\cx/\cy/\rot in {1/-\off/\off/\rotone, 2/\off/\off/\rottwo,
                               3/-\off/-\off/\rotthree, 4/\off/-\off/\rotfour} {
    \begin{scope}[shift={(\cx,\cy)}, rotate=\rot]
      \shade[inner color=black!10, outer color=black!3]
        (0,0) ellipse [x radius=\rx, y radius=\ry];
      \draw[part] (0,0) ellipse [x radius=\rx, y radius=\ry];
      \foreach \i in {1,...,\n} {
        \pgfmathsetmacro\vy{\hgt*(2*(\i-1)/(\n-1)-1)}
        \pgfmathsetmacro\vx{\rx*sqrt(1-(\vy/\ry)^2)}
        \coordinate (v\q-\i) at (\vx,\vy);
      }
      \node at ({-\rx-0.42},0) {$X_{\q}$};
    \end{scope}
  }
    \foreach \i/\j/\col in {
      5/1/ethree, 6/3/ethree, 4/1/efour, 5/5/efour, 1/2/ethree, 4/4/ethree, 6/5/ethree, 4/5/ethree, 6/1/ethree,
      3/1/ethree, 1/3/efour, 2/2/ethree, 4/2/ethree, 1/6/ethree, 3/5/ethree, 2/6/ethree, 6/4/ethree, 2/5/ethree,
      5/2/ethree, 3/4/ethree, 1/5/ethree, 4/3/efour, 2/3/ethree, 1/1/ethree, 2/1/efour, 5/3/ethree, 2/4/ethree,
      5/4/ethree, 3/6/ethree, 4/6/ethree, 3/2/efour, 1/4/efour, 3/3/ethree, 6/2/ethree, 5/6/efour, 6/6/efour}
      \draw[\col, cross] (v1-\i) -- (v4-\j);
    \foreach \i/\j/\col in {
      2/1/efour, 4/2/ethree, 6/6/ethree, 4/1/ethree, 6/4/efour, 6/2/ethree, 4/6/efour, 2/5/ethree, 1/1/ethree,
      3/2/ethree, 1/5/ethree, 5/4/efour, 1/6/ethree, 2/6/efour, 6/1/ethree, 6/3/ethree, 5/5/efour, 3/3/ethree,
      3/1/efour, 3/6/ethree, 4/4/ethree, 6/5/ethree, 5/3/ethree, 2/4/ethree, 4/3/ethree, 5/2/ethree, 1/4/ethree,
      2/2/ethree, 5/6/ethree, 2/3/ethree, 1/2/efour, 1/3/efour, 3/4/ethree, 4/5/ethree, 3/5/ethree, 5/1/ethree}
      \draw[\col, cross] (v2-\i) -- (v3-\j);
    \foreach \i/\j/\col in {
      1/1/eone, 2/2/eone, 2/4/efour, 5/2/eone, 1/6/efour, 2/1/eone, 4/1/eone, 2/6/eone, 2/3/eone,
      6/5/eone, 6/4/eone, 6/6/eone, 6/1/eone, 3/3/eone, 3/5/eone, 2/5/efour, 3/2/efour, 4/4/eone,
      3/1/eone, 5/3/eone, 1/4/eone, 5/5/eone, 1/3/efour, 3/6/eone, 3/4/eone, 1/2/eone, 4/2/eone,
      4/3/eone, 6/2/eone, 5/6/eone, 1/5/eone, 5/1/efour, 4/5/efour, 5/4/efour, 4/6/eone, 6/3/efour}
      \draw[\col, cross] (v1-\i) -- (v2-\j);
    \foreach \i/\j/\col in {
      5/6/eone, 6/4/efour, 3/1/eone, 2/3/eone, 6/3/eone, 1/3/eone, 5/5/eone, 3/6/eone, 2/4/eone,
      4/1/eone, 2/2/eone, 5/4/efour, 3/5/efour, 1/4/eone, 6/6/eone, 4/3/eone, 4/4/eone, 4/5/eone,
      4/6/efour, 3/4/eone, 6/1/eone, 4/2/eone, 6/5/efour, 5/2/efour, 2/6/eone, 1/1/eone, 6/2/eone,
      1/5/eone, 5/1/eone, 3/2/eone, 2/5/eone, 1/2/eone, 1/6/efour, 3/3/efour, 5/3/eone, 2/1/efour}
      \draw[\col, cross] (v3-\i) -- (v4-\j);
    \foreach \i/\j/\col in {
      1/6/etwo, 5/5/etwo, 4/3/efour, 2/5/etwo, 2/3/etwo, 5/2/etwo, 2/1/etwo, 4/1/efour, 6/1/etwo,
      3/5/efour, 1/2/efour, 1/1/etwo, 1/4/etwo, 2/2/efour, 1/5/etwo, 6/3/etwo, 3/4/efour, 2/6/etwo,
      5/6/etwo, 4/2/etwo, 1/3/etwo, 2/4/etwo, 3/2/etwo, 4/5/etwo, 6/2/etwo, 6/6/efour, 5/3/efour,
      4/4/etwo, 6/5/etwo, 5/4/etwo, 5/1/etwo, 3/6/etwo, 3/3/etwo, 4/6/etwo, 6/4/efour, 3/1/etwo}
      \draw[\col, cross] (v1-\i) -- (v3-\j);
    \foreach \i/\j/\col in {
      1/4/efour, 6/4/etwo, 2/1/etwo, 6/3/efour, 4/2/efour, 2/3/etwo, 1/1/efour, 6/6/etwo, 4/5/etwo,
      3/5/efour, 2/6/etwo, 4/6/etwo, 1/5/etwo, 3/6/efour, 3/4/etwo, 6/2/etwo, 3/1/etwo, 5/6/efour,
      6/5/etwo, 4/4/etwo, 1/6/etwo, 5/5/etwo, 5/4/etwo, 4/3/etwo, 4/1/etwo, 1/2/etwo, 5/3/etwo,
      5/1/etwo, 6/1/efour, 1/3/etwo, 2/4/etwo, 2/5/etwo, 5/2/etwo, 3/2/etwo, 2/2/efour, 3/3/etwo}
      \draw[\col, cross] (v2-\i) -- (v4-\j);
    \foreach \i/\j/\col in {
      1/6/etwo, 1/5/ethree, 2/6/eone, 1/4/etwo, 3/6/etwo, 2/5/ethree, 3/5/eone, 1/3/eone, 2/4/etwo,
      4/6/ethree, 2/3/eone, 3/4/ethree, 5/6/etwo, 1/2/ethree, 4/5/eone}
      \draw[\col] (v1-\i) to[out=135, in=135, looseness=\loose] (v1-\j);
    \foreach \i/\j/\col in {
      1/6/eone, 1/5/ethree, 2/6/etwo, 3/6/ethree, 1/4/etwo, 2/5/etwo, 3/5/etwo, 2/4/ethree, 1/3/eone,
      4/6/eone, 3/4/ethree, 1/2/eone, 5/6/eone, 4/5/etwo, 2/3/ethree}
      \draw[\col] (v2-\i) to[out=45, in=45, looseness=\loose] (v2-\j);
    \foreach \i/\j/\col in {
      1/6/etwo, 2/6/ethree, 1/5/eone, 2/5/ethree, 3/6/etwo, 1/4/eone, 4/6/ethree, 2/4/eone, 3/5/etwo,
      1/3/ethree, 2/3/ethree, 5/6/eone, 3/4/eone, 1/2/etwo, 4/5/etwo}
      \draw[\col] (v3-\i) to[out=225, in=225, looseness=\loose] (v3-\j);
    \foreach \i/\j/\col in {
      1/6/etwo, 2/6/eone, 1/5/ethree, 3/6/ethree, 2/5/eone, 1/4/etwo, 4/6/eone, 2/4/ethree, 1/3/etwo,
      3/5/eone, 1/2/eone, 5/6/etwo, 2/3/etwo, 3/4/ethree, 4/5/ethree}
      \draw[\col] (v4-\i) to[out=315, in=315, looseness=\loose] (v4-\j);
  \foreach \q in {1,...,4} \foreach \i in {1,...,\n} \node[vertex] at (v\q-\i) {};
\end{tikzpicture}
\endgroup